\documentclass[11pt]{amsart}

\numberwithin{equation}{section}

\begin{document}

\baselineskip=13pt

\title[]{On the variety of Lagrangian subalgebras, II}
\author[S. Evens and J.-H. Lu]{Sam Evens and Jiang-Hua Lu}

\address{Department of Mathematics, University of Notre Dame, Notre 
Dame, 46556
\newline \indent Department of Mathematics, University of Hong Kong, Pokfulam Road, Hong Kong}

\email{evens.1@nd.edu, jhlu@maths.hku.hk}

\begin{abstract}

\noindent VERSION FRAN{\c{C}}AISE: Motiv{\'e} par le th{\'e}or{\`e}me
 de Drinfeld
sur les espaces de Poisson homog{\`e}nes, nous {\'e}tudions la vari{\'e}t{\'e}
${\mathcal L}$ 
des sous-alg{\`e}bres de Lie Lagrangiennes de 
${\mathfrak g}\oplus{\mathfrak g}$
pour ${\mathfrak g}$,
une alg{\`e}bre de Lie complexe semisimple.
Soit $G$ le groupe adjointe de ${\mathfrak g}$.
Nous montrons que les adh{\'e}rences des $(G\times G)$-orbites dans 
${\mathcal L}$ sont les vari{\'e}t{\'e}s sph{\'e}riques et lisses.
Aussi, nous classifions les composantes irr{\'e}ductibles de ${\mathcal L}$
et nous montrons qu'elles sont lisses.
Nous employons quelques m{\'e}thodes de M. Yakimov pour donner
une nouvelle description et une 
nouvelle preuve de la classification de Karolinsky
des orbites diagonales de $G$ dans ${\mathcal L}$, quel, comme
cas sp{\'e}cial, donne la classification de Belavin-Drinfeld
des ${\mathfrak r}$-matrices quasitriangulaires de ${\mathfrak g}$.
En outre, ${\mathcal L}$ poss{\`e}de une structure de Poisson
canonique, et nous calculons son rang {\`a} chaque point and
nous d{\'e}crivons sa d{\'e}composition en feuilles symplectiques
en termes des intersections des orbites des deux sous-groupes de
$G\times G$.\vskip 5mm

\noindent ENGLISH VERSION: 
Motivated by Drinfeld's theorem on Poisson homogeneous spaces, we
study the variety ${\mathcal L}$ of Lagrangian subalgebras of ${\mathfrak g} \oplus {\mathfrak g}$
for a complex semi-simple Lie algebra ${\mathfrak g}$.
 Let $G$ be the adjoint group of ${\mathfrak g}$. We show that the
$(G \times G)$-orbit closures in ${\mathcal L}$ are smooth spherical varieties.
We also classify the irreducible components of ${\mathcal L}$ and show that
they are smooth. Using some methods of M. Yakimov,
we give a new description and proof of Karolinsky's classification of
the diagonal $G$-orbits in ${\mathcal L}$, which, as a special case,
recovers the Belavin-Drinfeld classification of
quasi-triangular r-matrices on ${\mathfrak g}$.
Furthermore, ${\mathcal L}$ has a canonical Poisson structure, and we compute its
rank at each point and describe its symplectic leaf decomposition
in terms of intersections of orbits of two subgroups of $G\times G$.

 \end{abstract}

\maketitle

\tableofcontents

\newtheorem{thm}{Theorem}[section]
\newtheorem{lem}[thm]{Lemma}
\newtheorem{prop}[thm]{Proposition}
\newtheorem{cor}[thm]{Corollary}
\newtheorem{rem}[thm]{Remark}
\newtheorem{exam}[thm]{Example}
\newtheorem{exam-nota}[thm]{Example-Notation}
\newtheorem{nota}[thm]{Notation}
\newtheorem{dfn}[thm]{Definition}
\newtheorem{ques}[thm]{Question}
\newtheorem{eq}{thm}
\newtheorem{dfn-nota}[thm]{Definition-Notation}
\newtheorem{claim}[thm]{Claim}
\newtheorem{dfn-lem}[thm]{Lemma-Definition}

\newcommand{\rw}{\rightarrow}
\newcommand{\lrw}{\longrightarrow}
\newcommand{\rhu}{\rightharpoonup}
\newcommand{\lhu}{\leftharpoonup}
\newcommand{\Map}{\longmapsto}
\renewcommand{\qed}{\begin{flushright} {\bf Q.E.D.}\ \ \ \ \
                  \end{flushright} }
\newcommand{\beqa}{\begin{eqnarray*}}
\newcommand{\eeqa}{\end{eqnarray*}}
\newcommand{\dia}{\begin{flushright}  $\diamondsuit$ \ \ \ \ \
                  \end{flushright} }

\newcommand{\ZR}{Z_{\Bbb R}}
\newcommand{\fls}{\fl_{d,\sigma}}
\newcommand{\flstriv}{\fl_\sigma}
\newcommand{\flphi}{\fl_\phi}
\newcommand{\Lphi}{L_\phi}
\newcommand{\flsigma}{\fl_\sigma}
\newcommand{\cphi}{c\circ \phi}
\newcommand{\Proj}{\rm Proj}

\newcommand{\flsC}{\fl_{d,\sigma, {\Bbb C}}}
\newcommand{\flsnull}{{\fl}_{d,{\sigma}_0}}
\newcommand{\flsnullC}{{\fl}_{d,{\sigma}_0, {\Bbb C}}}
\newcommand{\frds}{{\fr}_{d,\sigma}}
\newcommand{\Rds}{{R}_{d,\sigma}}
\newcommand{\os}{{\cal O}_{\Bbb C}}
\newcommand{\cO}{{\cal O}}
\newcommand{\ctilde}{{\tilde{c}}}
\newcommand{\wtilde}{{\tilde{w}}}

\newcommand{\la}{\mbox{$\langle$}}
\newcommand{\ra}{\mbox{$\rangle$}}
\newcommand{\ot}{\mbox{$\otimes$}}

\newcommand{\id}{\mbox{${\rm id}$}}
\newcommand{\Fun}{\mbox{${\rm Fun}$}}
\newcommand{\End}{\mbox{${\rm End}$}}
\newcommand{\Hom}{\mbox{${\rm Hom}$}}
\newcommand{\Ker}{\mbox{${\rm Ker}$}}
\renewcommand{\Im}{\mbox{${\rm Im}$}}

\newcommand{\Xa}{\mbox{$X_{\alpha}$}}
\newcommand{\Ya}{\mbox{$Y_{\alpha}$}}

\newcommand{\fc}{\mbox{${\mathfrak c}$}}
\newcommand{\fv}{\mbox{${\mathfrak v}$}}
\newcommand{\fd}{\mbox{${\mathfrak d}$}}
\newcommand{\fe}{\mbox{${\mathfrak e}$}}
\newcommand{\fa}{\mbox{${\mathfrak a}$}}
\newcommand{\ft}{\mbox{${\mathfrak t}$}}
\newcommand{\fk}{\mbox{${\mathfrak k}$}}
\newcommand{\fg}{\mbox{${\mathfrak g}$}}
\newcommand{\fq}{\mbox{${\mathfrak q}$}}
\newcommand{\fl}{\mbox{${\mathfrak l}$}}
\newcommand{\fs}{\mbox{${\mathfrak s}$}}
\newcommand{\fsl}{\mbox{${\fs\fl}$}}
\newcommand{\flp}{\mbox{${\mathfrak l}_{p}$}}
\newcommand{\fh}{\mbox{${\mathfrak h}$}}
\newcommand{\fn}{\mbox{${\mathfrak n}$}}
\newcommand{\fo}{\mbox{${\mathfrak o}$}}
\newcommand{\fp}{\mbox{${\mathfrak p}$}}
\newcommand{\fr}{\mbox{${\mathfrak r}$}}
\newcommand{\fb}{\mbox{${\mathfrak b}$}}
\newcommand{\fz}{\mbox{${\mathfrak z}$}}
\newcommand{\fm}{\mbox{${\mathfrak m}$}}
\newcommand{\fu}{\mbox{${\mathfrak u}$}}
\newcommand{\fbp}{\mbox{${\mathfrak b}_{+}$}}

\newcommand{\fgDelta}{\mbox{$\fg_{\Delta}$}}
\newcommand{\fnSminus}{\mbox{${\fn}_{S}^{-}$}}
\newcommand{\fnplus}{\mbox{${\fn}^+$}}
\newcommand{\fnminus}{\mbox{${\fn}^-$}}
\newcommand{\fnsigmaplus}{\mbox{${\fn_\sigma}$}}
\newcommand{\fnsigmaminus}{\mbox{${\fn_\sigma}_-$}}
\newcommand{\fnDeltaplus}{\mbox{${\fn_\Delta}_+$}}
\newcommand{\fnDeltaminus}{\mbox{${\fn_\Delta}_-$}}
\newcommand{\fmsDelta}{\mbox{${{\fm}_\Delta}_s$}}
\newcommand{\Geta}{\mbox{$G_\eta$}}
\newcommand{\Peta}{\mbox{$P_\eta$}}
\newcommand{\Gmin}{\mbox{${G^{-\tau_d}}$}}
\newcommand{\GDelta}{\mbox{$G_\Delta$}}
\newcommand{\ftdiag}{\mbox{${\ft}_{\Delta}$}}
\newcommand{\fmprime}{\mbox{${\fm}^\prime$}}
\newcommand{\fhsplit}{\mbox{${\fh}_{s}$}}
\newcommand{\Phiplushs}{\mbox{${\Sigma}^{+}(\fhsplit)$}}
\newcommand{\Phiplus}{\mbox{${\Sigma}^{+}$}}

\newcommand{\hs}{\hspace{.2in}}
\newcommand{\hhs}{\hspace{.15in}}

\newcommand{\fsk}{\mbox{${\fs}_{\frak k}$}}
\newcommand{\Ar}{\mbox{${\rm Aut}_{\frak r}$}}
\newcommand{\Ag}{\mbox{${\rm Aut}_{\frak g}$}}
\newcommand{\AS}{\mbox{${\rm Aut}_{\fmSss}$}}
\newcommand{\Ints}{\mbox{${\rm Int}_{\frak g}$}}
\newcommand{\Outr}{\mbox{${\rm Aut}_{D({\frak r})}$}}
\newcommand{\Outg}{\mbox{${\rm Aut}_{D({\frak g})}$}}
\newcommand{\daut}{\mbox{$d$}}
\newcommand{\dgamma}{\mbox{${\gamma}_d$}}
\newcommand{\donegamma}{\mbox{${\gamma}_{d_1}$}}
\newcommand{\dSTgamma}{\mbox{${\gamma}_{d}$}}
\newcommand{\dtheta}{\mbox{$\tau_d$}}
\newcommand{\dDelta}{\mbox{${\Delta,d}$}}
\newcommand{\Zdoneaut}{\mbox{${Z_{d_1}}$}}
\newcommand{\Zdaut}{\mbox{${Z_{d}}$}}
\newcommand{\ZdR}{Z_{d, {\Bbb R}}}
\newcommand{\Zdeta}{Z_{d, \eta}}
\newcommand{\Zdetao}{Z_{d, \eta_1}}

\newcommand{\fpS}{\mbox{$\fp_{ S}$}}
\newcommand{\fpSprime}{\mbox{${\fp}_{ S}^\prime$}}
\newcommand{\fmS}{\mbox{${\fm}_{S}$}}
\newcommand{\fmSprime}{\mbox{${\fm}_{ S}^\prime$}}
\newcommand{\fnS}{\mbox{${\fn}_{ S}$}}
\newcommand{\fpT}{\mbox{$\fp_{T}$}}
\newcommand{\fmT}{\mbox{${\fm}_{ T}$}}
\newcommand{\fnT}{\mbox{${\fn}_{T}$}}
\newcommand{\fnSone}{\mbox{${\fn}_{ S_1}$}}
\newcommand{\fnStwo}{\mbox{${\fn}_{ S_2}$}}
\newcommand{\fnTone}{\mbox{${\fn}_{ T_1}$}}
\newcommand{\fnTtwo}{\mbox{${\fn}_{ T_2}$}}
\newcommand{\fzS}{\mbox{${\fz}_{S}$}}
\newcommand{\fzT}{\mbox{${\fz}_{ T}$}}
\newcommand{\fzSone}{\mbox{${\fz}_{{ S}_1}$}}
\newcommand{\fzone}{\mbox{${\fz}_1$}}
\newcommand{\fmss}{\mbox{$\overline{{\fm}}$}}
\newcommand{\fmSss}{\fg_S}
\newcommand{\fmTss}{\fg_T}

\newcommand{\fmSv}{\fm_{S_w}}
\newcommand{\MSv}{G_{S_w}}
\newcommand{\fmSssv}{\fg_{S_w}}

\newcommand{\fmSssprime}{\mbox{$\overline{{\fm}_{S^\prime}}$}}
\newcommand{\Mss}{\mbox{$\overline{M}$}}
\newcommand{\fmssone}{\mbox{$\overline{{\fm_1}}$}}
\newcommand{\fnprime}{\mbox{${\fn}^\prime$}}
\newcommand{\fpprime}{\mbox{${\fp}^\prime$}}
\newcommand{\fzprime}{\mbox{${\fz}^\prime$}}
\newcommand{\fmssprime}{\mbox{${\fmss}^\prime$}}

\newcommand{\fmSigss}{\mbox{${\fm}_{\sigma, 1}$}}
\newcommand{\fgreal}{\mbox{${\fg}_0$}}
\newcommand{\PS}{\mbox{$P_{ S}$}}
\newcommand{\NS}{\mbox{$N_{ S}$}}
\newcommand{\MS}{\mbox{$M_{ S}$}}
\newcommand{\MSS}{G_S}
\newcommand{\MTS}{G_T}

\newcommand{\cR}{\cal R}
\newcommand{\Meta}{\mbox{$M_\eta$}}
\newcommand{\Netaplus}{\mbox{${N_\eta}+$}}
\newcommand{\Netaminus}{\mbox{${N_\eta}_-$}}
\newcommand{\Cleta}{\mbox{${{\C}^l}_{\eta}$}}
\newcommand{\Rleta}{\mbox{${{\R}^l}_{\eta}$}}
\newcommand{\Rl}{\mbox{${{\R}^l}$}}
\newcommand{\cstar}{\mbox{${\C}^*$}}
\newcommand{\fmsdeltaeta}{\mbox{$(\fmsDelta)^\eta$}}
\newcommand{\Lagr}{\mbox{${\mathcal L}$}}
\newcommand{\LinLagr}{{{\mathcal L}_{\rm space}}}
\newcommand{\Lagrfd}{\Lagr(\fd)}
\newcommand{\Gr}{{\rm Gr}(n, \fd)}
\newcommand{\RGr}{{\rm Gr}(n,\fg)}
\newcommand{\CGr}{{\rm Gr}(n,\fg \oplus \fg)}

\newcommand{\LagrfzS}{\mbox{${\mathcal L}_{{\mathfrak z}_{ S}}$}}
\newcommand{\LagrfzSC}{\mbox{${\mathcal L}_{{\mathfrak z}_{{ S},
{\Bbb C}}}$}}
\newcommand{\LagrfzSplus}{\mbox{$
{\mathcal L}_{{\mathfrak z}_{ S,1}}$}}
\newcommand{\LagrfzSminus}{\mbox{$
{\mathcal L}_{{\mathfrak z}_{ S,-1}}$}}
\newcommand{\LagrfzSeps}{\mbox{$
{\mathcal L}_{{\mathfrak z}_{ S,\epsilon}}$}}
\newcommand{\zeLagr}{\mbox{${\Lagr}_0$}}
\newcommand{\kLagr}{\mbox{${\Lagr}_{{\mathfrak k}}$}}
\newcommand{\lag}{\Lagr}
\newcommand{\kag}{\kLagr}
\newcommand{\supp}{{\rm supp}}
\newcommand{\fmst}{\fm_{\sigma}^{\theta_{d,\sigma}}}

\newcommand{\dw}{\dot{w}}
\newcommand{\du}{\dot{u}}
\newcommand{\dv}{\dot{v}}
\newcommand{\eps}{\epsilon}
\newcommand{\Vcplx}{\mbox{$V_0$}}
\newcommand{\ea}{\mbox{$E_{\alpha}$}}
\newcommand{\eb}{\mbox{$E_{-\alpha}$}}

\newcommand{\C}{\mbox{${\mathbb C}$}}
\newcommand{\Z}{\mbox{${\bf Z}$}}
\newcommand{\killing}{\ll \, , \, \gg}

\newcommand{\R}{\mbox{${\mathbb R}$}}
\renewcommand{\a}{\mbox{$\alpha$}}

\newcommand{\lara}{\la \, , \, \ra}
\newcommand{\llgg}{\ll \, , \, \gg}
\newcommand{\Ad}{{\rm Ad}}
\newcommand{\cF}{{\mathcal F}(S,T, d)}
\newcommand{\cFe}{{\mathcal F}^\epsilon(S, T, d)}
\renewcommand{\sl}{{\mathfrak sl}}

\newcommand{\CW}{{}^C\!W}
\newcommand{\flst}{\fl_{{\rm st}}}
\newcommand{\Pist}{\Pi_{{\rm st}}}
\newcommand{\Rone}{R^{\dot{w}}_{S_w, T_w, wd}}
\newcommand{\rone}{{\mathfrak r}^{\dot{w}}_{S_w, T_w, wd}}

\newcommand{\PSTw}{P_S \cap (\dw P_{T}^{-} \dot{w}^{-1})}

\newcommand{\swt}{{}^S\!W^T}
\newcommand{\Adw}{{\rm Ad}_{{\dot w}}}

\newcommand{\Grm}{{\rm Gr}(m, \fmSss \oplus \fmTss)}
\newcommand{\Grms}{{\rm Gr}(m, \fmSss \oplus \fmSss)}
\newcommand{\lagez}{\lag_{{\rm space}}^{\epsilon}(\fz_S \oplus \fz_T)}
\newcommand{\lagz}{\LinLagr(\fz_S \oplus \fz_T)}
\newcommand{\PSPT}{P_S \times P_{T}^{-}}
\newcommand{\ZdM}{Z_d(\MSS)}
\newcommand{\GPGP}{G/P_S \times G/P_{T}^{-}}
\newcommand{\mm}{\fmSss \oplus \fmTss}

\newcommand{\RSTd}{R_{S, T, d}}
\newcommand{\flstdv}{\fl_{S, T, d, V}}
\newcommand{\Advd}{\Ad_{\dot{v}} \gamma_d}
\newcommand{\fnsvd}{\fn_{S(v,d)}}
\newcommand{\Nsvd}{N_{S(v,d)}}
\newcommand{\MSvd}{M_{S(v,d)}}
\newcommand{\fzsvd}{\fz_{S(v,d)}}
\newcommand{\fgsvd}{\fg_{S(v,d)}}
\newcommand{\fmsvd}{\fm_{S(v,d)}}
\newcommand{\fnsvdp}{\fn_{S(v,d)}^{-}}
\newcommand{\Adv}{\Ad_{\dot{v}}}
\newcommand{\tpsi}{\tilde{\psi}}

\section{Introduction}
\label{sec_intro}

Let $\fd$ be a $2n$-dimensional Lie algebra over $k = \R$ or $\C$, 
together with a symmetric,  non-degenerate, 
and ad-invariant bilinear form $\la \,, \, \ra$.  A Lie subalgebra $\fl$ of $\fd$
is said to be {\it Lagrangian} if $\fl$ is maximal isotropic with
respect to $\la \, , \, \ra$, i.e., if
$\dim_k \fl = n$ and if $\la x, y \ra = 0$ for all
$x, y \in \fl$. By a {\it Lagrangian splitting} of $\fd$ 
we mean a direct sum
decomposition $\fd = \fl_1 + \fl_2$, where $\fl_1$ and $\fl_2$ are
two Lagrangian subalgebras of $\fd$. Denote by ${{\mathcal L}}(\fd)$ the set
of all Lagrangian subalgebras of $\fd$. It is an algebraic subvariety of the
Grassmannian ${\rm Gr}(n, \fd)$ of $n$-dimensional subspaces of $\fd$, and 
every connected 
Lie group $D$ with Lie algebra $\fd$ acts on ${{\mathcal L}}(\fd)$ via
the adjoint action of $D$ on $\fd$. We
proved in
\cite{e-l:reallag} that each Lagrangian
splitting $\fd = \fl_1 + \fl_2$ gives rise to a Poisson structure
$\Pi_{{\mathfrak l}_1, {\mathfrak l}_2}$ on ${{\mathcal L}}(\fd)$, making
${{\mathcal L}}(\fd)$ into a Poisson variety. Moreover, if $L_1$ and $L_2$ are
the connected subgroups of $D$ with Lie algebras
$\fl_1$ and $\fl_2$ respectively, all the $L_1$ and $L_2$-orbits in
${{\mathcal L}}(\fd)$ are Poisson submanifolds of $\Pi_{{\mathfrak l}_1,
{\mathfrak l}_2}$.

The above construction in \cite{e-l:reallag} was motivated by the work
of Drinfeld \cite{dr:homog} on Poisson homogeneous spaces. Indeed, a Lagrangian
splitting $\fd = \fl_1 + \fl_2$ of $\fd$ gives rise to the Manin triple
$(\fd, \fl_1, \fl_2)$, which in turn defines
Poisson structures $\pi_1$ and $\pi_2$ on $L_1$ and $L_2$
respectively, making them into Poisson Lie groups \cite{k-s:quantum}.
A Poisson space $(M, \pi)$ is said to be $(L_1, \pi_1)$-homogeneous
if $L_1$ acts on $M$ transitively and if the action map
$L_1 \times M \to M$ is Poisson. In \cite{dr:homog}, Drinfeld
constructed an $L_1$-equivariant map $M \to {{\mathcal L}}(\fd)$
for every $(L_1, \pi_1)$-homogeneous
Poisson space $(M, \pi)$ and proved that $(L_1, \pi_1)$-homogeneous
Poisson spaces correspond to $L_1$-orbits in ${{\mathcal L}}(\fd)$
in this way. The Poisson
structure $\Pi_{{\mathfrak l}_1, {\mathfrak l}_2}$ on ${{\mathcal L}}(\fd)$
is constructed in such a way that the Drinfeld map $M \to {{\mathcal L}}(\fd)$
is a Poisson map. In many cases, the Drinfeld map $M \to {{\mathcal L}}(\fd)$
is a local diffeomorphism onto its image.
Thus we can think of $L_1$-orbits in ${{\mathcal L}}(\fd)$ as
as models for $(L_1, \pi_1)$-homogeneous Poisson spaces. For this reason, it is
 interesting to study the geometry of the variety ${{\mathcal L}}(\fd)$,
the $L_1$ and $L_2$-orbits in ${{\mathcal L}}(\fd)$, and the Poisson 
structures $\Pi_{{\mathfrak l}_1, {\mathfrak l}_2}$ on 
${{\mathcal L}}(\fd)$.

There are many examples of Lie algebras $\fd$ with symmetric, non-degenerate, and
ad-invariant bilinear forms. The geometry of 
${{\mathcal L}}(\fd)$
is different from case to case. Moreover, there can be many 
Lagrangian splittings
for a given $\fd$, resulting in many Poisson structures on
${{\mathcal L}}(\fd)$.

\begin{exam}
\label{exam_fg-real}
{\em
Let $\fg$ be a complex semi-simple Lie algebra 
and regard it as a real Lie algebra, and
let $\la \, , \, \ra$ be the imaginary
part of the Killing form of $\fg$. The geometry of ${{\mathcal L}}(\fg)$ 
was studied in
\cite{e-l:reallag}. In particular, we determined the irreducible components of
${{\mathcal L}}(\fg)$ and classified the $G$-orbits in ${{\mathcal L}}(\fg)$, 
where $G$ is the adjoint group of $\fg$. Let $\fg = \fk + \fa + \fn$ be
an Iwasawa decomposition
of $\fg$. Then both $\fk$ and $\fa + \fn$ are Lagrangian,
so $\fg = \fk + (\fa + \fn)$ is a Lagrangian splitting,
resulting in a Poisson
structure on ${{\mathcal L}}(\fg)$ which we denote by $\pi_0$. 
Many interesting Poisson manifolds appear as $G$ or $K$-orbits 
inside ${{\mathcal L}}(\fg)$,
where $K$ is the connected subgroup of $G$ with Lie algebra $\fk$. 
Among such Poisson manifolds are the flag manifolds of $G$ and the compact
symmetric spaces associated to real forms of $G$. Detailed studies of the
Poisson geometry of these Poisson structures and some applications to Lie theory
have been given in \cite{lu:coor}, \cite{lu:cdyb}, \cite{e-l:harm}, and \cite{f-l:kk0}.
For example, a flag manifold $X$ of $G$ can be identified with a certain
$K$-orbit in ${{\mathcal L}}(\fg)$. The resulting Poisson structure $\pi_0$ 
on $X$ is
called the {\it Bruhat-Poisson structure} because its symplectic
leaves are Bruhat cells in $X$. In \cite{lu:coor} and \cite{e-l:harm},
we established connections between the Poisson geometry of $\pi_0$
and the harmonic forms on $X$ constructed by Kostant
\cite{ko:63} in 1963, and we gave a Poisson geometric interpretation of
the Kostant-Kumar approach \cite{k-k:integral} to Schubert  calculus
on $X$.
}
\end{exam}

\begin{exam}
\label{exam_semi}
{\em
Let $\fg$ be any $n$-dimensional Lie algebra, and let
$\fd$ be the semi-direct product of $\fg$ and its dual space $\fg^*$. The
symmetric bilinear form
$\la x + \xi, \, y + \eta \ra = \la x, \eta \ra + \la y, \xi \ra$ for 
$x, y \in \fg$ and $\xi, \eta \in \fg^*$
is non-degenerate and ad-invariant. When $\fg$ is semi-simple,
Lagrangian subalgebras of $\fd$ are
not easy to classify (except for low dimensional cases), for,  as a sub-problem,
one needs to classify all abelian subalgebras of $\fg$. See \cite{k-s:quantum},
\cite{karolinsky-stolin},
and their references for more detail. The description
of the geometry of $\lag(\fd)$ in this case is an open problem.
}
\end{exam}

\smallskip

In this paper, we consider the complexification of Example
\ref{exam_fg-real}. Namely, let $\fg$ be a
complex semi-simple Lie algebra and $\fd = \fg \oplus \fg$ the direct
sum Lie algebra with the
bilinear form
\[
\la (x_1,x_2),(y_1,y_2)\ra =\ll x_1,y_1\gg- \ll x_2,y_2\gg, \hs x_1, 
x_2, y_1, y_2 \in \fg,
\]
where $\ll \, ,\, \gg$ is a fixed symmetric, non-degenerate,
and ad-invariant
bilinear form on $\fg$. The variety of Lagrangian subalgebras of $\fd$ with
respect to $\la \, , \, \ra$ 
will be denoted by $\lag$.

The classification of Lagrangian subalgebras of $\fd$   has been
given by  Karolinsky
\cite{karo:homog-diag}, 
and Lagrangian splittings of $\fg \oplus \fg$ have been classified 
 by Delorme \cite{delorme:manin-triples}. 
In this paper, we establish the first few steps in the 
study of the Poisson structures on $\lag$ defined
by Lagrangian splittings of $\fg \oplus \fg$.
Namely,  we will first describe the geometry of 
$\lag$ in the following terms:

1)  the $(G \times G)$-orbits in $\lag$ 
and their closures, where $G$ is the adjoint group of $\fg$;

2)  the irreducible components of $\lag$;

\noindent
We then look at the Poisson structure $\Pi_0$ on $\lag$ 
defined by the so-called {\it standard} Lagrangian splitting
  $\fd =\fg_\Delta + \fg_{{\rm st}}^{*}$, where $\fg_{\Delta} = 
\{(x,x): x \in \fg\}$
and $\fg_{\rm st}^*
\subset \fb \oplus \fb^-$, where $\fb$ and $\fb^-$ are
two opposite Borel subalgebras of $\fg$. Let 
$G_\Delta = \{(g, g): g \in G\}$, and let $B$ and $B^-$ be the Borel 
subgroups of $G$ with Lie algebras $\fb$ and $\fb^-$ respectively. 
We will compute the rank of $\Pi_0$ in $\lag$ and study
the symplectic leaf decomposition of $\lag$ with respect 
to $\Pi_0$ in terms of 
the intersections of $G_\Delta$ and $(B \times B^-)$-orbits in
$\lag$.

We regard the Poisson structure $\Pi_0$ as the most important
Poisson structure on $\lag$, especially since it is more closely
related to applications to other areas of mathematics.
The study of the 
symplectic leaf decomposition of the Poisson structure 
$\Pi_{{\mathfrak l}_1, {\mathfrak l}_2}$
defined by an arbitrary Lagrangian splitting $\fg \oplus \fg = \fl_1 + \fl_2$
will be carried out in \cite{lu-milen2} using results of 
\cite{lu-milen1}. This analysis is much more technically involved,
so we think it is worthwhile to present the important special
case of $\Pi_0$ separately. In particular, in computing the 
rank of $\Pi_0$ in $\lag$, we need to classify $G_\Delta$-orbits in $\lag$
and compute the normalizer subalgebra of $\fg_\Delta$
for every $\fl \in \lag$. The classification of $G_\Delta$-orbits in
$\lag$ follows directly from a special case of Theorem 2.2 in \cite{lu-milen1}.
However, the inductive procedure
used in the proof of Theorem 2.2 in \cite{lu-milen1} 
is rather involved. To make the paper self-contained and 
more comprehensible, we present a shorter inductive proof for
the classification of 
$G_\Delta$-orbits in $\lag$.
Our proof and the one in 
\cite{lu-milen1} are both adapted from
\cite{yakimov:leaves}. The description of the normalizer subalgebra
of $\fg_\Delta$ at $\fl \in \lag$ can be derived from 
Theorem 2.5 in \cite{lu-milen1}. Again for the purpose of
completeness and comprehensibility, we give a simple and  more direct 
computation of these subalgebras.

We point out that E. Karolinsky has  
in \cite{karo:homog-diag} given
a classification of $G_\Delta$-orbits in $\lag$ in different terms.
Our classification is more in line with that
of Lagrangian splittings in \cite{delorme:manin-triples},
and in particular, the Belavin-Drinfeld theorem \cite{b-d:factorizable}
on Lagrangian splittings of the form $\fg \oplus \fg = \fg_\Delta
+ \fl$ follows easily
from our classification, but does not seem to be easily derived
from Karolinsky's classification. 
A special case of our classification is
the classification of $G_\Delta$-orbits on the
wonderful compactification of $G$,
and this special case was first determined by
Lusztig in \cite{lusz1} and \cite{lusz2}, where it is used to
develop a theory of character sheaves for the compactifications.
Our arguments are somewhat different
from those in \cite{lusz1} and \cite{lusz2}.
It would be very interesting to
find connections between Poisson geometry and character sheaves.

We now give more details of the main results in this paper:

Following O. Schiffmann \cite{schiff:classcdyb}, we define a 
{\it generalized Belavin-Drinfeld triple}
(generalized BD-triple) to be a triple
$(S, T, d)$, where $S$ and $T$ are two subsets of the
set $\Gamma$ of vertices of the Dynkin diagram of $\fg$, and $d: S \to T$ is an isometry with respect to $\ll \, , \, \gg$.
For a generalized BD-triple $(S, T, d)$, let $P_S$ and 
$P_{T}^{-}$ be respectively the standard parabolic
subgroups of $G$ of type $S$ and opposite type $T$ (see Notation \ref{nota_main})
and Levi decompositions
$P_S=M_SN_S$ and $P_{T}^{-}=M_TN_{T}^{-}$. Let
$G_S$ and $G_T$ be the quotients of  $M_S$ and $M_{T}$ 
by their centers respectively, and let
$\chi_S: M_S \to G_S$ and $\chi_T: M_T \to G_T$ be the natural projections.
Let $\gamma_d: G_S \to G_T$ be  the group isomorphism induced by $d$, and 
define the subgroup $R_{S, T, d}$ of  $P_S \times P_{T}^{-}$ by
\[
R_{S, T, d} = \{(m_S, m_T) \in M_S \times M_T: \; \gamma_d(\chi_S(m_S))
=\chi_T(m_T)\} (N_S \times N_{T}^{-}).
\]
In Section \ref{sec_L}, we establish the following facts on $(G \times G)$-orbits
and their closures in $\lag$ (Proposition \ref{prop_GG-orbit-types}, Corollary 
\ref{cor_GG-orbit-spherical}, and Proposition \ref{prop_GG-orbits-closure}):

\vspace{.03in}
{\it 
1) Every  $(G \times G)$-orbit in $\lag$ is isomorphic to 
$(G \times G)/\RSTd$ for a generalized BD-triple $(S, T, d)$,
so the $(G \times G)$-orbit
types in $\lag$ correspond bijectively to generalized
BD-triples for $G$; Every  $(G \times G)$-orbit in $\lag$  
is a $(G \times G)$-spherical homogeneous space.
 
2) For a
generalized BD-triple $(S, T, d)$, the  closure   of a 
   $(G \times G)$-orbit of type $(S, T, d)$
is a fiber bundle over the flag manifold 
$G/P_S \times G/P_{T}^{-}$
whose fiber is isomorphic to a De Concini-Procesi compactification
of $G_S$.  
}

\vspace{.03in}
We also study in Section \ref{sec_L} the irreducible components
of $\lag$. We prove (Corollary \ref{cor_2-components}, Theorem
\ref{thm_varclosure}, and Theorem \ref{thm_irredcomp}):

\vspace{.03in}
{\it
1) The irreducible components of $\lag$ are roughly 
(see Theorem \ref{thm_irredcomp} for detail)
labeled by quadruples $(S, T, d, \eps)$, where
$(S, T, d)$ are generalized BD-triples  and
$\eps \in \{0, 1\}$;

2) The irreducible component corresponding to 
$(S, T, d, \eps)$ is a fiber  bundle over 
the flag manifold $G/P_S \times G/P_{T}^{-}$ with fiber isomorphic to 
the product of a De Concini-Procesi compactification of $G_S$ and
a homogeneous space of a special orthogonal group. In particular,
all the irreducible components of $\lag$ are smooth;

3) $\lag$ has two connected components.
}
\vspace{.03in}
 
In Section \ref{sec_classification-G-orbits}, we give the classification of
$G_\Delta$-orbits in $\lag$, and we compute the normalizer subgalgebra of
$\fg_\Delta$ at every $\fl \in \lag$. In Section \ref{sec_Pi0},
we compute the rank of  the 
Poisson structure $\Pi_0$. In particular,
we prove
(Theorem \ref{thm_rank-Pi0} and Theorem \ref{thm_H-translates-leaves}):

\vspace{.03in}
{\it
Every non-empty intersection
  of a 
$G_\Delta$-orbit ${{\mathcal O}}$ and a  $(B \times B^-)$-orbit
${{\mathcal O}}^\prime$ in $\lag$ is 
a regular Poisson variety with respect to
$\Pi_0$;
The Cartan subgroup $H_\Delta=\{(h, h): h \in H\}$ of $G_\Delta$, 
where $H = B \cap B^-$,
acts transitively on the set of symplectic leaves in
${{\mathcal O}} \cap {{\mathcal O}}^\prime$.
}
\vspace{.03in}

An interesting Poisson subvariety of $(\lag, \Pi_0)$ is the 
De Concini-Procesi compactification $Z_1(G)$ of $G$ as the closure of
the $(G \times G)$-orbit through $\fg_\Delta \in \lag$. 
Conjugacy classes in $G$ and their closures in $Z_1(G)$ are all 
Poisson subvarieties of $(Z_1(G), \Pi_{0})$. In particular, 
$\Pi_0$ restricted to a conjugacy class $C$ in $G$
is non-degenerate precisely on the intersection of $C$
with the open Bruhat cell $B^-B$ (see Corollary \ref{cor_conjclassrank}).
It will be particularly interesting to compare the 
Poisson structure $\Pi_0$ on the unipotent variety in $G$ with
the Kirillov-Kostant structure on the nilpotent cone in $\fg^*$.
Another interesting Poisson subvariety of $(\lag, \Pi_0)$ is 
the De Concini-Procesi compactification $X_\sigma$ 
of a complex symmetric space
$G/G^\sigma$ for an involutive automorphism $\sigma$ of $G$ 
(Proposition \ref{prop_orbitequalsdp}). 
Intersections of $G_\Delta$-orbits and $(B \times B^-)$-orbits
inside the  closed $(G \times G)$-orbits in $\lag$ 
are related to double Bruhat cells in $G$ (see Example \ref{exam_double-Bruhat}), 
and Kogan and Zelevinsky \cite{kogan-zelevinsky}
have constructed toric charts on some of the symplectic leaves
in these closed orbits.
It would be interesting to 
see how their methods can be applied to other symplectic leaves of $\Pi_0$.

{\bf Acknowledgment:} We would like to thank Milen Yakimov and 
Eugene Karolinsky for pointing out  errors in
a preliminary version of the paper. Discussions
with Milen Yakimov enabled us to improve earlier results and 
solve problems in more complete forms. We would also like to
thank Michel Brion, William Graham, and George McNinch for useful
comments and the referee for suggesting many corrections
and improvements. The second author is grateful to the 
Hong Kong University of Science and Technology for
its hospitality.
The first author was partially supported by (USA)NSF
grant DMS-9970102 and the second author by 
(USA)NSF grant DMS-0105195, HKRGC grants 701603 and 703304, and 
the New Staff Seeding Fund
at HKU. 
 
\section{The variety $\lag$ of Lagrangian subalgebras of $\fg\oplus\fg$}
\label{sec_L}
\smallskip

Throughout this paper, $\fg$ will be a 
complex semi-simple Lie algebra, and $\ll \, , \, \gg$ will be a 
fixed symmetric and
non-degenerate
ad-invariant bilinear form on $\fg$.
Equip the direct product Lie algebra
$\fg \oplus \fg$ with the symmetric non-degenerate ad-invariant bilinear form
\begin{equation}
\label{eq_I}
\la (x_1,x_2),(y_1,y_2)\ra =\ll x_1,y_1\gg- \ll x_2,y_2\gg, \hs
x_1, x_2, y_1, y_2 \in \fg.
\end{equation}
A Lagrangian subalgebra of $\fg \oplus \fg$ is by definition
an $n$-dimensional complex Lie subalgebra of $\fg
\oplus \fg$ that is isotropic with respect to $\la \, , \, \ra$.

\begin{nota}
\label{nota_L}
{\em We will denote by ${{\mathcal L}}$ the
variety of all Lagrangian subalgebras of $\fg \oplus \fg$
and by $\LinLagr(\fg \oplus \fg)$ the
variety of all $n$-dimensional isotropic {\it subspaces}
of $\fg \oplus \fg$.}
\end{nota}

Let $G$ be the adjoint group of $\fg$.  
In this section,
we will use the classification
of Lagrangian subalgebras of $\fg \oplus \fg$ 
by Karolinsky \cite{karo:homog-diag} to study the
$(G \times G)$-orbits and their closures in $\lag$, and
we will
determine the irreducible components of ${{\mathcal L}}$.

\subsection{Lagrangian subspaces}
\label{sec_lagrsubspace}

Let $U$ be a finite-dimensional complex vector space with a
 symmetric and non-degenerate bilinear form $\la \, , \, \ra$. A subspace $V$ of $U$ is
said to be Lagrangian
if $V$ is maximal isotropic with respect to $\lara$.
If $\dim U = 2n$ or $2n+1$, Witt's theorem
says that the dimension of a Lagrangian subspace of $U$
is $n$. Denote by $\LinLagr(U)$ the set of Lagrangian subspaces of $U$. 
It is a closed 
algebraic subvariety of ${\rm Gr}(n, U)$, the Grassmannian of
$n$-dimensional subspaces of $U$.

\begin{prop}[\cite{ACGH}, pp. 102-103]
\label{prop_indlagr}
Let $\dim U = 2n$ (resp.  $2n +1$) with $n > 0$. Then
$\LinLagr(U)$ is a smooth algebraic subvariety of ${\rm Gr}(n, U)$ with
two (resp. one) connected components,
each of which is isomorphic to 
$SO(2n, \C)/P$
(resp. $SO(2n+1, \C)/P)$) where $P$ has Levi factor isomorphic to $GL(n, \C)$.
Moreover, $\LinLagr(U)$ has complex dimension
$\frac{n(n-1)}{2}$ (resp. $\frac{n(n+1)}{2}$).
When $\dim U = 2n$, two Lagrangian subspaces
$V_1$ and $V_2$ are in the same connected component of $\LinLagr(U)$
if and only if $\dim(V_1) - \dim(V_1\cap V_2)$ is even.
\end{prop}

\begin{nota}
\label{nota_L0-L1}
{\em For $U = \fg \oplus \fg$ with the bilinear form
$\lara$ in (\ref{eq_I}),  let $\lag^0$ be the intersection of
$\Lagr$ with the connected component of $\LinLagr(\fg \oplus \fg)$
containing the diagonal of $\fg \oplus \fg$. The intersection of $\lag$
with the other connected component of $\LinLagr(\fg \oplus \fg)$
will be denoted by $\lag^1$.
}
\end{nota}

Let $\fh$ be a Cartan subalgebra of $\fg$, and let $\fn$ and $\fn^-$ be the opposite
nilpotent subalgebras of $\fg$ corresponding to a  choice of positive roots for
$(\fg, \fh)$. For a Lagrangian subspace
$V$ of $\fh \oplus \fh$ with respect to $\lara$, let
$\fl_V = V + \{(x, y):  x \in \fn, y \in \fn^-\}.$
Then $\fl_V \in \lag$. It is
easy to see from Proposition \ref{prop_indlagr} that $\fl_{V_1}$
and $\fl_{V_2}$ are in the same connected component of
$\LinLagr(\fg \oplus \fg)$ if and only if
$V_1$ and $V_2$ are in the same connected component of
$\LinLagr(\fh \oplus \fh)$. In particular,
  $\Lagr^1$ is non-empty.

\subsection{Isometries}
\label{sec_isometries}

We collect some results on automorphisms that will be used in
later sections.

\begin{nota}
\label{nota_Cartan}
{\em Throughout this paper, we
fix a Cartan subalgebra $\fh$ and a choice $\Sigma^+$ of
positive roots in the set $\Sigma$ of all roots
of $\fg$ relative to $\fh$, and let $\fg = \fh + \sum_{\alpha \in
\Sigma} \fg_\alpha$ be the root decomposition.
Let $\Gamma$ be the
set of simple roots in $\Sigma^+$. For $\alpha \in \Sigma$,  let
$H_\alpha \in \fh$ be such that $\ll H_\alpha, H \gg = \alpha(H)$
for all $H \in \fh$. For $\alpha \in \Sigma^+$, fix
root vectors $\ea \in \fg_\alpha$ and $\eb \in \fg_{-\alpha}$
such that $\ll \ea, \eb \gg = 1$.  For a subset 
$S$ of $\Gamma$, let $[S]$ be the set of roots in the linear span of $S$, and let
$G_S$ be the adjoint group of the semisimple Lie algebra $\fg_S$ given by
\[
\fg_S ={\rm span}_{{\mathbb C}}\{H_\alpha: \alpha \in S\}  + \sum_{\alpha
\in [S]} \fg_\alpha \subset \fg.
\]
}
\end{nota}

For $S, T \subset \Gamma$, we are interested in Lie algebra
isomorphisms $\fmSss \to \fmTss$ that preserve the restrictions of
$\llgg$ to $\fmSss$ and $\fmTss$. We will simply
refer to this property as preserving $\llgg$. 

\begin{dfn}
\label{dfn_d-isometry}
{\em
Let $S, T \subset \Gamma$. By an isometry
from $S$ to $T$ we mean a bijection $d: S \to T$ such that
$\ll d\alpha, d\beta \gg = \ll \alpha, \beta \gg$ for
 all $\alpha, \beta \in S$,
where $\ll \alpha, \beta \gg = \ll H_\alpha, H_\beta \gg$.
Let $I(S,T)$ be the set of all isometries from $S$ to $T$.
Following \cite{schiff:classcdyb}, a
triple $(S, T, d)$, where $S, T \subset \Gamma$ and $d \in I(S, T)$,
will also be called a {\it generalized Belavin-Drinfeld (generalized 
BD-)triple for $G$}.
}
\end{dfn}

\begin{lem}
\label{lem_leviisomorphism} For any $S, T \subset \Gamma$ and 
$d \in I(S, T)$, there is a unique isomorphism
$\gamma_d: \fmSss \to \fmTss$  
such that
\begin{equation}
\label{eq_dSTgamma}
\dSTgamma(E_{\alpha})=E_{d(\alpha)}, \hs
\dSTgamma(H_{\alpha})=H_{d(\alpha)}, \hs \forall \alpha \in S.
\end{equation}
Moreover,  $\gamma_d$
preserves
$\ll \, ,\, \gg$, and for every Lie algebra
isomorphism $\mu:\fmSss\to\fmTss$ preserving $\llgg$, there is a unique
$d \in I(S, T)$ and a unique $g \in \MSS$ such that
$\mu = \gamma_d {\Ad}_g$.
\end{lem}

\noindent
{\bf Proof.}  Existence and uniqueness of $\gamma_d$ is by Theorem 2.108 
in \cite{knapp}.  For $\alpha \in \Sigma^+$, let $\lambda_\alpha, 
\mu_\alpha \in \C$ be such that
$\gamma_d(E_\alpha) = \lambda_\alpha E_{d\alpha}$ and
$\gamma_d(E_{-\alpha}) = \mu_\alpha E_{-d\alpha}$. By applying
$\gamma_d$ to the identity $[E_\alpha, E_{-\alpha}] = H_\alpha$ we get
$\lambda_\alpha \mu_\alpha = 1$ for every $\alpha \in \Sigma^+$.
It follows that  $\gamma_d$ preserves $\ll \,, \, \gg$.
Suppose that $\mu:\fmSss\to\fmTss$ is a Lie algebra isomorphism
preserving $\llgg$.
Let $d_1$ be any isomorphism  from the Dynkin diagram of
$\fmSss$ to the Dynkin diagram of $\fmTss$.
Let $\gamma_{d_1}: \fmSss \to \fmTss$ be defined as in (\ref{eq_dSTgamma}).
Then $\nu := \gamma_{d_1}^{-1}  \mu$ is an
automorphism of $\fmSss$. Recall that there is a short exact sequence
\[
1 \, \longrightarrow \MSS \, \longrightarrow \, \AS \,
{\longrightarrow} \, {\rm Aut}_S\, \longrightarrow \, 1,
\]
where $\AS$ is the group of 
automorphisms of $\fmSss$, and ${\rm Aut}_S$ is the group of 
automorphisms of the Dynkin diagram of $\fmSss$.
Let $d_2 \in {\rm Aut}_S$ be the image of $\nu$ under the map
$\AS \to {\rm Aut}_S$ and write
$\nu = \gamma_{d_2} {\rm Ad}_g$ for some $g\in \MSS$. Thus 
$\mu = \gamma_{d_1} \gamma_{d_2} \Ad_g = \gamma_{d_1d_2} \Ad_g$.
Since $\mu$ and $\Ad_g$ are isometries of $\llgg$,
so is $\gamma_{d_1 d_2}$. Thus $d := d_1 d_2 \in I(S, T)$ is an isometry, and
$\mu = \gamma_d \Ad_g$.

Uniqueness of $d$ and $g$ follows from the fact that if
$g_0\in \MSS$ preserves a Cartan subalgebra and acts as the
identity on all simple root spaces, then $g_0$ is the identity element.
\qed

\begin{dfn}
\label{dfn_type-d-1}
{\em A $\ll \, , \, \gg$-preserving Lie algebra isomorphism $\mu: \fmSss \to \fg_T$
is said to be of {\it type $d$} for $d \in I(S, T)$
if $d$ is the unique element in $I(S, T)$ such that
$\mu = \gamma_d \Ad_g$ for some $g \in \MSS$.
}
\end{dfn}

\subsection{Karolinsky's classification}
\label{sec_karoclass}

Karolinsky \cite{karo:homog-diag} has classified the Lagrangian
subalgebras of $\fg\oplus\fg$
with respect to the bilinear form $\la \, , \, \ra$ given in
(\ref{eq_I}). We recall his results now.

\begin{nota}
\label{nota_parabolic}
{\em
For a parabolic subalgebra  $\fp$ of $\fg$,
let $\fn$ be its nilradical and $\fm:=\fp/\fn$ its Levi factor.
Let $\fm=[\fm, \fm] + \fz$ be the decomposition of $\fm$ into the direct sum of
its derived
algebra $[\fm, \fm]$ and its center $\fz$. Recall that $[\fm, \fm]$ is 
semisimple and that the restrictions of the bilinear form $\ll \, , \, \gg$ to
$\fm$ and $\fz$ are both non-degenerate.
If $\fpprime$ is another parabolic subalgebra, denote its nilradical,
Levi factor,
and
center of Levi factor, etc. by $\fnprime$, $\fmprime$, and
$\fzprime$, etc.. 
When speaking of Lagrangian subspaces of $\fz\oplus\fzprime$, we
mean with respect to the restriction of $\la \, , \, \ra$ to $\fz\oplus\fzprime$.
}
\end{nota}

\begin{dfn}
\label{dfn_admquad} {\em  A quadruple $(\fp,\fpprime,\mu,V)$ is
called admissible if $\fp$ and $\fpprime$ are parabolic subalgebras of
$\fg$,
$\mu: [\fm, \fm]\to [\fm^\prime, \fm^\prime]$ is a Lie algebra
isomorphism preserving $\ll \, , \, \gg$, and $V$ is a Lagrangian
subspace of $\fz\oplus \fzprime$.}
\end{dfn}

If $(\fp,\fpprime,\mu,V)$ is admissible, set
\[
\fl(\fp,\fpprime,\mu,V):= \{(x, x^\prime): x \in \fp, x^\prime \in \fpprime,
\mu(x_{[{\mathfrak m}, {\mathfrak m}]} )= x^{\prime}_{[{\mathfrak m}^\prime, 
{\mathfrak m}^\prime]},
(x_{{\mathfrak z}}, x^{\prime}_{{\mathfrak z}^\prime}) \in V\} 
\subset \fg \oplus
\fg,
\]
where for $x \in \fp$, $x_{[{\mathfrak m}, {\mathfrak m}]}\in [\fm, \fm]$ and
$x_{{\mathfrak z}} \in \fz$ are respectively the $[\fm, \fm]$- and
$\fz$-components of $x + \fn \in \fp / \fn = [\fm, \fm] + \fz$. 
We use similar notation for $x^\prime \in \fpprime$.

\begin{thm}[\cite{karo:homog-diag}]
\label{thm_karo_cplx}
$\fl(\fp,\fpprime,\mu,V)$ is a Lagrangian subalgebra if
$(\fp,\fpprime,\mu,V)$ is admissible, and 
every Lagrangian subalgebra of $\fg\oplus\fg$ is of the form
$\fl(\fp,\fpprime,\mu,V)$
for some admissible quadruple.

\end{thm}

\subsection{Partition of $\Lagr$}
\label{sec_partition}

\begin{nota}
\label{nota_main}
{\em 
Recall the fixed choice $\Sigma^+$ of positive roots from Notation \ref{nota_Cartan}.
Set
\[
\fn = \sum_{\alpha \in \Sigma^+} \fg_\alpha, \hs
\fn^-= \sum_{\alpha \in \Sigma^+} \fg_{-\alpha}.
\]
A parabolic subalgebra $\fp$ of $\fg$
is said to be {\it standard} if it contains the Borel subalgebra
$\fb:=\fh +\fn.$
Recall also  that, for $S \subset \Gamma$, 
$[S]$ denotes
the set of roots in the linear span of $S$.  Let
\[
\fmS = \fh + {\sum}_{\alpha \in [S]}^{} {\fg}_\alpha, \hspace{.2in} 
\fnS = {\sum}_{\alpha\in {\Sigma}^+ -[S]}^{}
 {\fg}_\alpha, \hs \fn_{S}^{-} = \sum_{\alpha \in \Sigma^+ - [S]}
\fg_{-\alpha}
\]
and $\fpS = \fmS + \fnS, \, \fp_{S}^{-} = \fm_S + \fn_{S}^{-}.$
We will refer to $\fpS$ as the standard parabolic subalgebra of
$\fg$ defined by $S$ and $\fp_{S}^{-}$
the {\it opposite} of $\fp_S$. A parabolic subalgebra of $\fp$ of $\fg$
is said to be of {\it type $S$} if
$\fp$ is conjugate to $\fp_S$ and of the 
{\it opposite-type} $S$ if $\fp$ is conjugate to $\fp_{S}^{-}$.
Let $\fmSss = [\fmS, \fmS]$ as in Notation \ref{nota_Cartan} and
\[
\fh_S = \fh\cap\fmSss= {\rm span}_{{\mathbb C}} \{H_\alpha: \alpha \in
[S]\}, \, \hs \, \fz_S = \{x \in \fh: \, \alpha(x) = 0, \,
\forall \alpha \in S\}.
\]
Then
$\fp_S = \fz_S + \fmSss + \fn_S$ and $\fp_{S}^{-} = \fz_S + \fmSss + \fn_{S}^{-}$.
The connected subgroups of $G$ with Lie algebras $\fp_S, \fp_{S}^{-},
 \fm_{S},  \fn_S$ and $\fn_{S}^{-}$ will be respectively denoted by
$P_S, P_{S}^{-},  M_S, N_S$ and $N_{S}^{-}$. Correspondingly
we have $P_S = M_S N_S$ and $P_{S}^{-} = M_S N_{S}^{-}$,
where $M_S \cap N_S =\{e\} = M_S \cap N_{S}^{-}$.
Recall that $\MSS$ denotes the adjoint group of $\fmSss$. Let $\chi_S:
P_S \to G_S$ be the composition of the projection from $P_S = M_SN_S$ to $M_S$
along $N_S$ and the projection $M_S \to G_S$. 
The similarly defined projection 
from $P_{S}^{-}$ to $\MSS$ will also be denoted by $\chi_S$.
}
\end{nota}

Returning to the notation in Notation \ref{nota_parabolic}, we have

\begin{dfn-lem}
\label{lem_leviisomorphism-1}
Let $(\fp, \fpprime, \mu)$ be a triple, where $\fp$ and
$\fpprime$ are parabolic subalgebras of $\fg$, and
$\mu: [\fm, \fm] \to [\fm^\prime, \fm^\prime]$ is a Lie algebra isomorphism
preserving $\llgg$. Assume that $\fp$ is of type $S$ and $\fpprime$ is of
opposite-type
$T$. Let $g_1, g_2 \in G$ be such that $\Ad_{g_1} \fp = \fp_S$ and $\Ad_{g_2} 
\fp^\prime = \fp_{T}^{-}$. Let
$\overline{\Ad_{g_1}}$ and $\overline{\Ad_{g_2}}$ be the induced Lie
algebra isomorphisms
\[
\overline{\Ad_{g_1}}: \,  [\fm, \fm] \lrw \fmSss, \hs
\overline{\Ad_{g_2}}: \,   [\fm^\prime, \fm^\prime]\lrw \fmTss.
\]
If $\mu^\prime := \overline{\Ad_{g_2}}\circ \mu \circ 
(\overline{\Ad_{g_1}})^{-1}:
\fmSss \to \fmTss$ is of type $d \in I(S, T)$ as in Definition \ref{dfn_type-d-1}, 
we say that
$(\fp, \fpprime, \mu)$ is of {\it type} $(S, T, d)$.
The
type of $(\fp, \fpprime, \mu)$ is independent of the choice of $g_1$
and $g_2$.
\end{dfn-lem}

\noindent
{\bf Proof.} If $h_1$ and $h_2$ in $G$ are such that $\Ad_{h_1} \fp = \fp_S$
and $\Ad_{h_2} \fp^\prime = \fp_{T}^{-}$, then there exist $p_S \in P_S$ and
$p_{T}^{-} \in P_{T}^{-}$ such that $h_1 = p_S g_1$ and $h_2 = p_{T}^{-} g_2$. 
Thus
\[
\mu^{\prime \prime} := \overline{\Ad_{h_2}} \circ \mu
\circ (\overline{\Ad_{h_1}})^{-1}
= \overline{\Ad_{p_{T}^{-}}} \circ \mu^\prime \circ (\overline{\Ad_{p_S}})^{-1}.
\]
The action of $\overline{\Ad_{p_S}}$ on $\fmSss$ is by definition the
adjoint action of $\chi_S(p_S) \in \MSS$ on $\fmSss$.
Similarly for the action of $\overline{\Ad_{p_{T}^{-}}}$ on $\fmTss$. Thus
by Definition \ref{dfn_type-d-1}, the two maps $\mu^\prime$
and $\mu^{\prime \prime}$ have the same type.
\qed
   
We are now ready to partition $\lag$. Recall the definitions of
${{\mathcal L}}^0$
and ${{\mathcal L}}^1$ in Notation \ref{nota_L0-L1}.

\begin{dfn}
\label{dfn_normaut}
{\em
Let $S,T\subset \Gamma, d\in I(S,T)$, and $\eps \in \{0, 1\}$. Define
 $\Lagr^\epsilon (S,T,d)$ to be the set of all
Lagrangian subalgebras
$\fl(\fp,\fpprime,\mu,V)$ such that

1) $\fl(\fp,\fpprime,\mu,V) \in \lag^\epsilon$;

2) $(\fp, \fpprime, \mu)$ is of type $(S, T, d)$.

\noindent
We say that $\fl \in \Lagr$ is of {\it type} $(\eps, S, T, d)$ if $\fl
\in \Lagr^\epsilon (S,T,d)$.}
\end{dfn}

It is clear that 
\begin{equation}
\label{eq_disjoint-L}
\Lagr = \bigcup_{\epsilon \in \{0, 1\}} \, \,
\bigcup_{S,T \subset \Gamma,
d\in I(S,T)} \Lagr^\epsilon(S,T,d)
\end{equation}
is a disjoint union, 
and each $\Lagr^\epsilon(S,T,d)$ is invariant under $G \times G$.
Set
\[
\fn_S \oplus \fn_{T}^{-} = \{(x, y): x \in \fn_S, y \in \fn_{T}^{-}\}
\subset \fg \oplus \fg,
\]
and for each $V \in \LinLagr(\fz_S \oplus \fz_T),$ set
\begin{equation}
\label{eq_flSTdV}
\fl_{S, T, d, V} = V + (\fn_S \oplus \fn_{T}^{-}) +
\{(x, \gamma_d(x)): x \in \fmSss\} \in \Lagr.
\end{equation}
Note that $\dim \fz_S = \dim \fz_T$ because $\fh_S \cong \fh_T$ and
$\fh = \fz_S + \fh_S = \fz_T + \fh_T$ are direct sums.

\begin{dfn-lem}
\label{dfn_two-eps}
For  $V_1, V_2 \in \LinLagr(\fz_S \oplus \fz_T)$,  
$\fl_{S, T, d, V_1}$ and $\fl_{S, T, d, V_2}$
are in the same connected component of  $\LinLagr(\fg \oplus \fg)$
if and only if $V_1$ and $V_2$ are in the same connected component
of $\LinLagr(\fz_S \oplus \fz_T)$. For $\epsilon = \{0, 1\}$,
let
\[
\lag_{{\rm space}}^{\epsilon}(\fz_S \oplus \fz_T) = \{V \in
\LinLagr(\fz_S \oplus \fz_T): \, \fl_{S, T, d, V} \in \lag^\epsilon\}.
\]
\end{dfn-lem}

\noindent
{\bf Proof.} The statement follows from Proposition \ref{prop_indlagr} and
the fact
\[
\dim(\fl_{S,T,d,V_1}) - \dim(\fl_{S, T, d, V_1} \cap \fl_{S, T, d, V_2}) = 
\dim (V_1) - \dim (V_1 \cap V_2).
\]
\qed

\begin{prop}
\label{prop_Le-union}
For any generalized BD-triple $(S, T,d)$
 and $\eps \in \{0, 1\}$, 
\begin{equation}
\label{eq_Leunion}
\lag^{\epsilon} (S, T, d) = \bigcup_{V \in \lag_{{\rm space}}^{\epsilon}
(\fz_S \oplus \fz_T)}
(G \times G) \cdot \fl_{S, T, d, V} \hs \hs (\mbox{disjoint union}).
\end{equation}
\end{prop}

\noindent
{\bf Proof.} By Definition
\ref{dfn_normaut}, every $(G \times G)$-orbit in
$\lag^{\epsilon} (S, T, d)$ passes through an
$\fl_{S, T, d, V}$ for some $V \in
\lag_{{\rm space}}^{\epsilon}(\fz_S \oplus \fz_T)$. If
$V_1, V \in \lag_{{\rm space}}^{\epsilon}(\fz_S \oplus \fz_T)$
are such that $\fl_{S, T, d, V_1}= \Ad_{(g_1, g_2)}
\fl_{S, T, d, V}$, then $(g_1, g_2)$ normalizes
$\fn_S \oplus \fn_{T}^{-}$, so $(g_1, g_2) \in P_S \times P_{T}^{-}$,
and it follows that $V_1 = V$.
\qed

\subsection{$(G \times G)$-orbits in $\lag$}
\label{sec_GG-orbits}

The following theorem follows immediately from Proposition
\ref{prop_Le-union}
and the decomposition of $\lag$ in (\ref{eq_disjoint-L}).

\begin{thm}
\label{thm_GG-orbits-1}
 Every $(G \times G)$-orbit in $\lag$ passes through
an $\fl_{S, T, d, V}$ for a unique
 quadruple $(S, T, d, V)$, where $S, T \subset \Gamma, d \in I(S, T)$ 
 and
$V \in \LinLagr(\fz_S \oplus \fz_T)$.
\end{thm}

For each $S, T \subset \Gamma$ and $d \in I(S, T)$, define the group
(see Notation  \ref{nota_main})
\begin{equation}
\label{eq_RSTd}
R_{S, T, d} := \{(p_S, p_{T}^{-}) \in P_S \times P_{T}^{-}: \, \,
\gamma_d(\chi_S(p_S))=\chi_T(p_{T}^{-}) \} \subset \PSPT.
\end{equation}

\begin{prop}
\label{prop_GG-orbit-types}
Let $S, T \subset \Gamma, d \in I(S, T)$, and
$V \in \LinLagr(\fz_S \oplus \fz_T)$.

1) The
$(G \times G)$-orbit in $\lag$ through $\fl_{S, T, d, V}$
is isomorphic to $(G \times G)/R_{S, T, d}$ and it has dimension
$n - z$, where $n = \dim \fg$ and $z = \dim \fz_S$.

2) $(G \times G)\cdot \fl_{S, T, d, V}$ fibers over
$G/P_S \times G/P_{T}^{-}$ with fiber isomorphic to $G_S$.
\end{prop}

\noindent
{\bf Proof.} It is routine to check that the stabilizer of
$\fl_{S, T, d, V}$ is $R_{S, T, d}$, and the dimensional
formula follows. The fiber in 2) may be identified with
$(P_S \times P_{T}^-)/R_{S, T, d}$, which may be identified
with $G_S$ via the map 
\[
(P_S \times P_{T}^-)/R_{S, T, d} \lrw G_S: \; \; (p_S, p_{T}^-)\mapsto
\gamma_{d}^{-1}(\chi_T(p_{T}^{-})) (\chi_S(p_S))^{-1}
\]
\qed

\begin{rem}
\label{rem_orbit-type}
{\em
It follows that the number of orbit types for $G \times G$ in
$\lag$ is exactly the number of generalized BD-triples for $G$.
}
\end{rem}

\begin{lem}
\label{lem_RSTd-connected}
$R_{S, T, d}$ is connected.
\end{lem}

\noindent
{\bf Proof.} The projection map $p:R_{S, T, d}\to P_S$, $(p_S,p_{T}^{-})
\mapsto
p_S$ is surjective and has fiber $N_{T}^{-} \times Z_T$, where $Z_T$
is the center of the group $M_T$. Clearly $N_{T}^{-}$
 is connected, and $Z_T$ is connected by Proposition
8.1.4 of \cite{carter:lietype}. Since $P_S$ is
connected, it follows that $R_{S, T, d}$ is also connected.\qed

\subsection{$(B \times B^{-})$-orbits in $\lag$}
\label{sec_BB-orbits}

Let $B$ and $B^-$ be the Borel subgroups of $G$ with
Lie algebras
$\fb = \fh + \fn$  and $\fb^- = \fh + \fn^-$ respectively.
By Proposition \ref{prop_GG-orbit-types}, 
to classify $(B \times B^-)$-orbits in $\lag$,
it suffices to consider
$(B \times B^{-})$-orbits in $(G \times G)/\RSTd$ for any  
generalized BD-triple
$(S, T, d)$.

Let $W$ be the Weyl group of $\Sigma$. For $S \subset \Gamma$,
let $W_S$ be the subgroup of $W$ generated by elements in $S$, and let 
$W^S$ be the set of minimal length
representatives of elements in the cosets in $W/W_S$.
Fix a representative $\dot{w}$ in $G$ for  each $w \in W$.
The following assertion can be proved in 
the same way as Lemma 1.3 in \cite{sp:IC}. 
It also follows directly from Proposition 8.1 of \cite{lu-milen1}.

\begin{prop}
\label{prop_BB-finite-orbits}
Let $(S, T, d)$ be an generalized BD-triple for $G$. Then
every  $(B \times B^-)$-orbit in $(G \times G)/\RSTd$ 
is of the form 
$(B \times B^-)(\dot{w}, \dot{v}) \RSTd$ for a unique
pair $(w, v) \in W \times W^T$.
\end{prop}

\begin{cor}
\label{cor_BB-finite-orbits}
Every $(B \times B^-)$-orbit in $\lag$ goes through
exactly one point in $\lag$ of the form
$\Ad_{(\dot{w}, \dot{v})} \flstdv$, where
$(S, T, d)$ is a generalized BD-triple,
$V \in \lag_{{\mathrm space}}(\fz_S \oplus \fz_T)$, and $(w, v) \in W \times W^T$.
\end{cor}

Since each $(G \times G)$-orbit in $\lag$ has finitely many 
$(B \times B^-)$-orbits,   at least one
of them is open. 

\begin{cor}
\label{cor_GG-orbit-spherical}
All  $(G \times G)$-orbits in $\lag$ are $(G \times G)$-spherical 
homogeneous spaces.
\end{cor}

\subsection{The De Concini-Procesi compactifications
$Z_d(G)$ of $G$}
\label{sec_wonderful-of-G}

In this section, we will consider the closure in
$\lag$ of some special $(G \times G)$-orbits.
Namely, when $S = T = \Gamma$ and $d \in
I(\Gamma, \Gamma)$, we have the graph $\fl_{\gamma_d}$ of
$\gamma_d: \fg \to \fg$ as a point in $\lag$:
\begin{equation}
\label{eq_graph-gammad}
{\fl}_{\gamma_d} = \{(x, \gamma_d(x)): \, x \in \fg\}.
\end{equation}
The $(G \times G)$-orbit in $\lag$ through $\fl_{\gamma_d}$
can be identified with $G$ by the map
\begin{equation}
\label{eq_embedding-Gd}
I_d: \, G \lrw (G \times G) \cdot \fl_{\gamma_d}: \, \,
g \Map \{(x, \, \gamma_d \Ad_g (x)): \, x \in \fg\}.
\end{equation}
The identification $I_d$ is $(G \times G)$-equivariant if we equip
$G$ with the action of $G \times G$ given by
\begin{equation}
\label{eq_GG-on-G}
(g_1, g_2) \cdot g = \gamma_{d}^{-1} (g_2) g g_{1}^{-1}.
\end{equation}

Since an orbit of an algebraic group on a variety is locally closed
(Section 8.3 in \cite{hum:linear-algebraic-groups}),
the orbit $(G \times G) \cdot {\fl}_{\gamma_d}$ has the same
closure in the Zariski topology and in the classical topology.
The closure
$\overline{(G \times G) \cdot {\fl}_{\gamma_d}}$, called a
{\it De Concini-Procesi compactification of $G$}, is a
smooth projective variety of dimension $n= \dim G$ (see
\cite[$\S 6$]{dp:compactification}). We denote this closure
by $Z_d(G)$.

It is known  \cite{dp:compactification} that $G\times G$ has
finitely many orbits in $\Zdaut(G)$ indexed by subsets of
 $\Gamma$. Indeed, for each
$S \subset \Gamma$, let   $\fl_{S, d}
\in \Lagr$ be given by
\begin{equation}
\label{eq_flSd}
\fl_{S,d}= {\fnS} \oplus \fn_{d(S)}^{-}
+ \{(x,\dgamma(x)): \, x\in \fmS \}.
\end{equation}
Choose $\lambda \in \fh$ such that $\alpha(\lambda) = 0$
for $\alpha \in S$ and $\alpha(\lambda) > 0$ for
$\alpha \in \Gamma - S$, and let 
$e^\lambda:\cstar \to H$ be the 
one parameter
subgroup of $H$ corresponding to $\lambda$.
Then it is easy to see that
\[
\lim_{t \to +\infty} \Ad_{(e^\lambda(t), e)} \fl_{\gamma_d} = \fl_{S, d}
 \in {\rm Gr}(n,
\fg \oplus \fg).
\]
Thus $\fl_{S,d} \in \Zdaut(G)$. It is easy to see that
$  \fl_{\gamma_d} \in \lag^\epsilon(\Gamma, \Gamma, d)$
for $\eps = ({\dim \fh - \dim \fh^{\gamma_d}})\mod 2$. Thus
$\fl_{S, d} \in \lag^\epsilon(S,d(S),d|_S)$ for the same value of $\eps$.

\begin{thm}
\label{thm_GGorbits}
 \cite{dp:compactification}
 For every $d \in I(\Gamma, \Gamma)$, $\Zdaut(G) =  \bigcup_{S \subset \Gamma}
 (G \times G) \cdot \fl_{S, d}$.
\end{thm}

\subsection{Closures of $(G \times G)$-orbits in $\lag$}
\label{sec_closures-GG-orbits}

Let $(S, T, d, V)$ be a quadruple where $(S, T, d)$ is a
generalized BD-triple and
$V \in \LinLagr(\fz_S \oplus \fz_T)$.
For $\fl_{S, T, d, V}$ given in (\ref{eq_flSTdV}),
we will now study the closure of the $(G \times G)$-orbit through
$\fl_{S, T, d, V}$ in ${\rm Gr}(n, \fg \oplus \fg)$. To this end,
let
${\rm Gr}(m, \fmSss \oplus \fmTss)$ be the Grassmannian
of $m$-dimensional
subspaces in $\fmSss \oplus \fmTss$, where
 $m = \dim \fmSss$.
For the Lie algebra isomorphism $\gamma_d: \fmSss \to \fmTss$
given in (\ref{eq_dSTgamma}), let
\[
\fl_{\gamma_d} = \{(x, \, \gamma_d(x)): \, x \in \fmSss\}
\in  {\rm Gr}(m, \fmSss \oplus \fmTss).
\]

\begin{nota}
\label{nota_ZdM}
{\em Let $Z_d(\MSS)$ be the closure of
 $(\MSS \times \MTS) \cdot \fl_{\gamma_d}$
in ${\rm Gr}(m, \fmSss \oplus \fmTss)$. Under the identification
\begin{equation}
\label{eq_embedMSS}
\MSS \lrw (\MSS \times \MTS) \cdot \fl_{\gamma_d}: \,\, g  \Map \{(x, \gamma_d \Ad_{g } x): \,
x \in \fmSss\},
\end{equation}
$\ZdM$ can be regarded as a $(\PSPT)$-equivariant
compactification of $\MSS$, where $\PSPT$ act on
$\MSS$ by
\begin{equation}
\label{eq_action-PSPT-on-MS-1}
(p_S, p_{T}^{-}) \cdot g_S = \gamma_{d}^{-1}(\chi_T(p_{T}^{-})) \, g_S \,
(\chi_S(p_S))^{-1}, \hs (p_S, p_{T}^{-}) \in \PSPT, \, g_S \in \MSS
\end{equation}
(see Notation \ref{nota_main}), and on
$\Grm$ via the group homomorphism
$\chi_S \times \chi_T: \PSPT \to \MSS \times \MTS$.
}
\end{nota}

\begin{prop}
\label{prop_GG-orbits-closure}
For every generalized BD-triple
$(S, T, d)$ and every $V \in \LinLagr(\fz_S \oplus \fz_T)$,

1) the closure $\overline{(G \times G) \cdot \fl_{S, T, d, V}}$ in
${\rm Gr}(n, \fg \oplus \fg)$
is a smooth subvariety of ${\rm Gr}(n, \fg \oplus \fg)$ of
dimension $n -z$, where $n = \dim \fg$ and
$z = \dim \fz_S$, and the map
\[
{\bf a}: \, \,  (G \times G) \times_{(P_S \times P_{T}^{-})}
Z_d(\MSS) \lrw
\overline{(G \times G) \cdot \fl_{S, T, d, V}}: \, \, 
[(g_1, g_2), \fl] \Map \Ad_{(g_1, g_2)} (V +
(\fn_S \oplus \fn_{T}^{-}) + \fl)
\]
is a $(G \times G)$-equivariant isomorphism;

2) $\overline{(G \times G) \cdot \fl_{S, T, d, V}}$
is the finite disjoint union
\[
\overline{(G \times G) \cdot \fl_{S, T, d, V}} = \bigcup_{S_1 \subset S} 
(G \times G) \cdot
\fl_{S_1, d(S_1), d_1, V_1},
\]
where for $S_1 \subset S$, we set  $d_1 = d|_{S_1}$, and
$V_1 = V + \{(x, \gamma_d(x)):  \, x \in \fh_S \cap \fz_{S_1}\}
\subset \fz_{S_1}
\oplus \fz_{T_1}.$
 \end{prop}

\noindent
{\bf Proof.} 
 Since $G/P_S \times G/P_{T}^{-}$ is
complete,  the image of
${\bf a}$ is closed. 
Since ${\bf a}[(e,e),\fl_{\gamma_d}]=\fl_{S,T,d,V}$ and
${\bf a}$ is $(G\times G)$-equivariant, ${\bf a}$ has dense image.
 Hence ${\bf a}$ is onto. 2) follows easily from the fact
that ${\bf a}$ is onto
and the description of orbits in $Z_d(\MSS)$ in 
Theorem \ref{thm_GGorbits}. Further, by this description
of orbits, for $\fl \in Z_d(\MSS)$, 
the nilradical of $V + (\fn_S \oplus \fn_{T}^{-}) + \fl$
contains $\fn_S \oplus \fn_{T}^-$. It follows that
if ${\bf a}[(g_1,g_2),\fl] \in {\bf a}[(e,e),Z_d(\MSS)]$,
then $(g_1,g_2)$ is in the normalizer of the nilradical
of $V + (\fn_S \oplus \fn_{T}^{-}) + \fl$. This normalizer,
and hence also $(g_1,g_2)$, is contained in
$P_S \times P_T^-$. Now suppose 
${\bf a}[(g_1,g_2),\fl]={\bf a}[(x_1,x_2),\fl_1]$ for
$x_1, x_2 \in G$ and $\fl_1 \in Z_d(\MSS)$. Then
 ${\bf a}[(x_{1}^{-1}g_1,x_{2}^{-1}g_2),\fl]=
{\bf a}[(e,e),\fl_1]$, so by the above $(x_{1}^{-1}g_1,
 x_{2}^{-1}g_2) \in P_S \times P_T^-$.
 It follows easily
that ${\bf a}$ is injective.

To show that the differential ${\bf a}_*$ of ${\bf a}$
is injective everywhere, it suffices to show ${\bf a}_*$ is
injective at all points $[(e,e),\fl]$ by 
$(G\times G)$-equivariance.
For $X\in \fg\oplus \fg$, let $\xi_X$ be the vector
field on $(G \times G) \times_{(P_S \times P_{T}^{-})} Z_d(\MSS)$
tangent to the $(G\times G)$-action. Recall that the tangent
space at a plane $U$ to the Grassmannian ${\rm Gr}(n, Y)$
may be identified with ${\rm Hom}(U,Y/U)$. Since $Z_d(\MSS)
\subset {\rm Gr}(\dim(\fg_S),\fg_S + \fg_T)$, it follows
that the tangent
space to $(G \times G) \times_{(P_S \times P_{T}^{-})} Z_d(\MSS)$
at $[(e,e),\fl]$ is spanned by the set of 
$\xi_X$ for $X\in \fn_S^-\oplus \fn_T$
and the set of $\phi$ in ${\rm Hom}(\fl,(\fg_S \oplus \fg_T)/\fl)$
tangent to $Z_d(\MSS)$.
We identify the tangent space to
${\rm Gr}(n, \fg\oplus \fg)$ at ${\bf a}[(e,e),\fl]$ with
${\rm Hom}(V + \fn_S\oplus \fn_T^- + \fl, 
 (\fg \oplus \fg)/(V + \fn_S\oplus \fn_T^- + \fl))$.
Then ${\bf a}_*(\xi_X)=\eta_X$, where for 
$Y \in V + \fn_S\oplus \fn_T^-+\fl$,
$\eta_X(Y)= {\rm ad}_X (Y) \ {\rm mod } \ (V + \fn_S\oplus \fn_T^-+\fl)$.
Moreover,  ${\bf a}_*(\phi)=\tilde{\phi}$, where $\tilde{\phi}(x)=\phi(x)$
for $x\in \fl$, and $\tilde{\phi}(y)=0$ if $y \in V + \fn_S\oplus \fn_T^-$.
We further identify $(\fg \oplus \fg)/(V + \fn_S\oplus \fn_T^-+\fl)$
with  $(\fz_S \oplus \fz_T)/V
+ \fn_S^-\oplus \fn_T +(\fg_S \oplus \fg_T)/\fl$.
Now suppose that ${\bf a}_*(\xi_X + \phi)=0$. Since
${\rm ad}_X (\fl)\subset \fn_S^-\oplus \fn_T$ and $\tilde{\phi}(\fl)
\subset (\fg_S \oplus \fg_T)/\fl$, it follows that $\tilde{\phi} = 0$.
Hence, $\phi = 0$, and $\xi_X = 0$ since ${\bf a}_*$ is injective on
$\fn_S^-\oplus \fn_T$.
By Corollary 14.10 of \cite{harr:agbook},
${\bf a}$ is an isomorphism.
\qed

Consider now the case when $S$ and $T$ are the empty set $\emptyset$,
so $d = 1$. By Theorem \ref{thm_GG-orbits-1},
every $(G \times G)$-orbit in
$\lag^0(\emptyset, \emptyset, 1) \cup
\lag^1(\emptyset, \emptyset, 1)$
goes through a unique
Lagrangian subalgebra of the form
\begin{equation}
\label{eq_flV-1}
\fl_V = V + (\fn \oplus \fn^-),
\end{equation}
where $V \in \LinLagr(\fh \oplus \fh)$.
 The following fact
follows immediately from Proposition \ref{prop_GG-orbits-closure}.

\begin{cor}
\label{cor_closed-orbits}
For every $V \in  \LinLagr(\fh \oplus \fh)$, the
$(G \times G)$-orbit through $\fl_V$ is isomorphic to
$G/B \times G/B^-$. These are the only closed
$(G \times G)$-orbits in $\lag$.
\end{cor}

\begin{cor}
\label{cor_2-components}
$\lag$ has two
connected components.
\end{cor}

\noindent
{\bf Proof.} In Section
\ref{sec_lagrsubspace}, we observed that $\lag$ has at least two connected
components, namely $\lag^0$ and $\lag^1$. Since every orbit of an
algebraic group on a variety has a closed orbit
in its boundary (see Section 8.3 in \cite{hum:linear-algebraic-groups}),
every point in $\lag$ is in the same connected component as $\fl_V$ for
some $V \in  \LinLagr(\fh \oplus \fh)$.
Thus $\lag$ has at most two connected components.
\qed

\subsection{The geometry of the strata $\Lagr^\epsilon(S,T,d)$}
\label{sec_geometry-strata}

Fix a generalized BD-triple $(S, T, d)$ and
an $\eps \in \{0, 1\}$.
Let
$P_S \times P_{T}^{-}$ act on $G_S$ by
(\ref{eq_action-PSPT-on-MS-1}) and 
trivially on $\lag_{{\rm space}}^{\epsilon} (\fz_S \oplus \fz_T)$. 
Consider the associated bundle
$(G \times G) \times_{(P_S \times P_{T}^{-})}
(\MSS \times \lag_{{\rm space}}^{\epsilon}(\fz_S \oplus \fz_T))$
over
$G/P_S \times G/P_{T}^{-}$ and the map
\begin{eqnarray}
\label{eq_a-2}
{\bf a}: & & (G \times G) \times_{(P_S \times P_{T}^{-})}
(\MSS \times \lag_{{\rm space}}^{\epsilon} (\fz_S \oplus \fz_T)) \lrw
\Lagr^\epsilon(S,T,d)\\
& & [(g_1, g_2), (g, V)] \Map
\Ad_{(g_1, g_2)} \fl_{g, V},
\end{eqnarray}
where $\fl_{g, V} = V + (\fn_S \oplus \fn_{T}^{-}) +
\{(x, \gamma_d \Ad_g (x)): x \in \fg_S\}$ for $g \in G_S$.

\begin{prop}
\label{prop_varopenpart}
For every $S, T \subset \Gamma, d \in I(S, T)$, and $\eps \in
\{0, 1\}$,
$\Lagr^\epsilon(S,T,d)$ is a smooth connected subvariety of
${\rm Gr}(n, \fg \oplus \fg)$ of dimension
$n + \frac{z(z-3)}{2}$, where $n = \dim \fg$ and $z=\dim{\fzS}$,
and the map ${\bf a}$ in (\ref{eq_a-2}) is
a $(G \times G)$-equivariant
isomorphism.
\end{prop}

\noindent{\bf Proof.} Consider the $(G \times G)$-equivariant projection
\begin{equation}
\label{eq_J}
J: \, \, \Lagr^\epsilon(S,T,d) \lrw G/P_S \times G/P_{T}^{-}:
\, \, \fl(\fp, \fpprime, \mu, V) \Map (\fp, \fpprime).
\end{equation}
Let $\cFe$ be the fibre of $J$ over the point $(\fp_S, \fp_{T}^{-})
\in G/P_S \times G/P_{T}^{-}$. 
By Lemma 4, p. 26
of \cite{slod:simpsing},
 the map
\[
(G \times G) \times_{(P_S \times P_{T}^{-})} \cFe \lrw
\Lagr^\epsilon(S,T,d):\, \,
[(g_1, g_2), \fl] \Map \Ad_{(g_1, g_2)} \fl
\]
is a $(G \times G)$-equivariant isomorphism.  By Lemma
\ref{lem_leviisomorphism},
\[
\cFe  =  \{\fl_{g, V}: \, \,  g \in \MSS, \, \, V \in
\lag_{{\rm space}}^{\epsilon}(\fz_S \oplus \fz_T)\}.
\]
 The identification $G_S \times
\lag_{{\rm space}}^{\epsilon}(\fz_S \oplus \fz_T) \to \cFe:
(g, V) \mapsto \fl_{g, V}$ is $(P_S \times P_{T}^{-})$-equivariant.
It follows that ${\bf a}$ is a $(G\times G)$-equivariant
isomorphism.
The dimension claim 
follows from
Propositions \ref{prop_indlagr} and \ref{prop_GG-orbit-types}. 
Smoothness and connectedness of
$\lag^\epsilon(S, T, d)$ follow easily.
\qed

\subsection{The geometry of the closure of
$\overline{\lag^\epsilon(S, T, d)}$}
\label{sec_geometry-closures-strata}

For a generalized BD-triple $(S, T, d)$ and
an $\eps \in \{0, 1\}$, recall that
$\ZdM$ is an $(\PSPT)$-equivariant compactification of $G_S$. 
Let again $\PSPT$ act trivially on $\lag_{{\rm space}}^{\epsilon} (\fzS\oplus \fz_T)$.
We omit the proof of the following Theorem since it is similar to
the proof of Theorem \ref{prop_GG-orbits-closure}.

\begin{thm}
\label{thm_varclosure}
For every generalized BD-triple $(S, T, d)$
and $\eps \in \{0, 1\}$, the closure
$\overline{\Lagr^\epsilon(S,T,d)}$
is a smooth algebraic variety of dimension
$n + \frac{z(z-3)}{2}$, where $n =\dim(\fg)$, $z=\dim{\fzS}$, and 
\begin{eqnarray}
\label{eq_a}
{\bf a}: & &
(G \times G) \times_{(P_S \times P_{T}^{-})} (Z_{{d}}(\MSS) \times
\lag_{{\rm space}}^{\epsilon} (\fzS\oplus \fz_T))
 \lrw  \overline{\lag^\epsilon(S, T, d)}\\
& & \; \; [(g_1, g_2), (\fl, V)] \Map
\Ad_{(g_1, g_2)} (V + (\fn_S \oplus \fn_{T}^{-}) + \fl)
\end{eqnarray}
 is a $(G \times G)$-equivariant isomorphism.
\end{thm}

\begin{cor}
\label{cor_Lebar-union}
For every generalized BD-triple $(S, T, d)$
and $\eps \in \{0, 1\}$, 
\begin{equation}
\label{eq_Lebar-union}
\overline{\lag^{\epsilon}(S, T, d)}=\bigcup_{V \in
{{\mathcal L}}^{\epsilon}_{{\mathrm space}}({{\mathfrak z}}_S
\oplus {{\mathfrak z}}_T)} \; \; \bigcup_{S_1 \subset S}
(G \times G) \cdot \fl_{S_1, d(S_1), d, V_1(V, S_1)}
\end{equation}
is a disjoint union, where for $S_1 \subset S$ and
$V \in \lag_{{\rm space}}^{\epsilon}(\fz_S \oplus \fz_T)$,
\[
V_1(V, S_1) = V + \{(x, \gamma_d(x)): x \in \fh_S \cap \fz_{S_1}\}
\subset \fz_{S_1} \oplus \fz_{d(S_1)}.
\]
\end{cor}

\begin{rem}
\label{rem_two-topologies-same-closure}
{\em
1). Since $Z_{{d}}(\MSS)$ is also the closure in
the classical topology
of $(\MSS \times \MTS)\cdot \fl_{\gamma_d}$
inside ${\rm Gr}(m, \fmSss \oplus \fmTss)$, 
$\lag^\epsilon(S, T, d)$ also has the same closure in the two
topologies of ${\rm Gr}(n, \fg \oplus \fg)$.

2). Since $Z_{{d}}(\MSS) - \MSS$ has dimension strictly lower than
$m = \dim \MSS$, it follows from the proof of Theorem
\ref{thm_varclosure} that $\overline{\lag^{\epsilon}(S, T, d)}
-\lag^\epsilon(S, T, d)$ is of strictly lower
dimension than $\overline{\lag^{\epsilon}(S, T, d)}$.
}
\end{rem}

 \subsection{Irreducible components of $\Lagr$}
\label{sec_irred}

Since each $\overline{\lag^\epsilon(S, T, d)}$ is smooth and
connected, it is a
closed irreducible subvariety of $\lag$.
Since
\[
\Lagr = \bigcup_{\epsilon \in \{0, 1\}} \; \;
\bigcup_{S,T \subset \Gamma, d \in I(S,T)}
 \overline{\Lagr^\epsilon(S, T, d)}
\]
is a finite union, the irreducible components of $\Lagr$ are those
$\overline{\Lagr^\epsilon(S,T, d)}$ that are not properly contained
in some other such set.

\begin{thm}
\label{thm_irredcomp}
$\overline{\Lagr^\epsilon(S,T,d)}$ is an irreducible component
of $\Lagr$ unless $|\Gamma -S| = 1$, $T =d_1(S)$ for some
$d_1 \in I(\Gamma, \Gamma)$, $d = d_1|_S$,
and
$\eps= (\dim\fh- \dim{\fh}^{\gamma_{d_1}}) \mod 2$.
\end{thm}

\noindent{\bf Proof.}
 When $(S, T, d, \eps)$ are as described in the proposition,
$\dim \fz_S = 1$, so 
$\Lagr^\epsilon(S,T,d)$ consists of a single
$(G \times G)$-orbit which lies in $Z_{d_1}(G)$ by
Theorem \ref{thm_GGorbits}.
We need to show that
this is the only nontrivial
case when the closure $\overline{\Lagr^\epsilon(S,T,d)}$
is contained in another $\overline{\Lagr^\epsilon(S_1, T_1, d_1)}$.

Assume that ${\Lagr^\epsilon(S,T,d)}$
 is in the boundary
of $\overline{\Lagr^\epsilon(S_1,T_1,d_1)}$.
Then by Corollary \ref{cor_Lebar-union},
$S\subset S_1$ and $T\subset T_1$.
By Remark  \ref{rem_two-topologies-same-closure},
$\dim {\Lagr^\epsilon(S,T,d)} < \dim \Lagr^\epsilon(S_1, T_1, d_1)$,
and thus
\[
\frac{1}{2}\dim(\fzS) (\dim(\fzS)-3) <
\frac{1}{2}\dim(\fz_{S_1}) (\dim(\fz_{S_1})-3)
\]
by  the dimension formula in Proposition \ref{prop_varopenpart}.
Since $S\subset S_1$, so $\dim(\fzS) \ge \dim(\fz_{S_1})$,
these two inequalities
imply that $\dim(\fz_{S_1})=0$ and $\dim(\fzS)=1$ or $2$.
In particular,
$S_1= T_1=\Gamma$,
so $\eps = (\dim \fh - \dim \fh^{\gamma_{d_1}})\mod 2$, and
$\overline{\Lagr^\epsilon(S_1,T_1,d_1)}= Z_{d_1}(G)$.

If $\dim(\fzS)=2$, $\Lagr^\epsilon(S,T,d)$ contains infinitely many
$(G \times G)$-orbits by Theorem \ref{thm_GG-orbits-1} and
Proposition \ref{prop_indlagr}. Since
$Z_{d_1}(G)$ has only finitely many $(G \times G)$-orbits,
$\Lagr^\epsilon(S,T,d)$ can not be contained in
$Z_{d_1}(G)$.
Assume that $\dim(\fzS)=1$. Then by Proposition
\ref{prop_varopenpart},
$\Lagr^\epsilon(S,T,d)$ is a single
$(G \times G)$-orbit. By the description of the
$(G \times G)$-orbits in $Z_{d_1}(G)$ in
Theorem \ref{thm_GGorbits}, 
$T$ and $d$ must be as described in the proposition.
\qed

\begin{exam}
\label{exam_sl2}
{\em 
For $\fg = \sl(2, \C)$, $\lag$ has two
irreducible components. One is
the De Concini-Procesi compactification $Z_{{\rm id}}(G)$ of $G = PSL(2, C)$
 which is isomorphic to $\C P^3$ (see
\cite{dp:compactification}), and the other is
 the closed $(G \times G)$-orbit
through $\fh_{\Delta} + (\fn  \oplus \fn^-)$, and is isomorphic to
$\C P^1 \times \C P^1$.

For $\fg = \sl(3, \C)$, there are four irreducible components
$Z_{{\rm id}}(G), Z_{d_1}(G), C_1$ and $C_2$, where
$Z_{{\rm id}}(G)$ and $Z_{d_1}(G)$ are the two
 De Concini-Procesi
compactifications of $G = PSL(3, \C)$ corresponding to the
identity and the non-trivial automorphism of the Dynkin
diagram of $\sl(3, \C)$, and $C_1$ and $C_2$
are the two components
$\lag_0(\emptyset, \emptyset, d)$ and $\lag_1(\emptyset, \emptyset, d)$.
Both $C_1$ and $C_2$ have dimension $7$. Moreover, $Z_{{\rm id}}(G) \cap C_1$
is a $6$-dimensional closed $(G \times G)$-orbit, and so is
$Z_{d_1}(G) \cap C_2$.
}
\end{exam}

\section{Classification of $G_\Delta$-orbits in $\lag$}
\label{sec_classification-G-orbits}

By Theorem \ref{thm_GG-orbits-1}, to describe
the $G_\Delta$-orbits in $\lag$, it suffices to describe
$(G_\Delta, \RSTd)$-double cosets in $G \times G$
for all generalized BD-triples
$(S, T, d)$ for $\Gamma$, where
$\RSTd$ is given by (\ref{eq_RSTd}). 
%In Section 
%\ref{sec_a-double-coset-theorem}, we will 
%define a class of subgroups $R$ of $G \times G$ that are
%slightly more general than the subgroups $\RSTd$, and 
%we will classify
%$(G_\Delta, R)$-double cosets in $G \times G$ for such
%an $R$.
A general double coset theorem in \cite{lu-milen1} 
classifies  $(R_{S^\prime, T^\prime, d^\prime}, 
R_{S, T, d})$-double
cosets in $G \times G$  for
two arbitrary generalized BD-triples $(S^\prime,T^\prime,d^\prime)$ and 
$(S, T, d)$. However, the proof in \cite{lu-milen1} is rather technically
 involved.
In this section, we present a simplified proof 
for the special case when
$R_{S^\prime, T^\prime, d^\prime} = G_\Delta$.
 Our method and the one used in 
\cite{lu-milen1} are both adapted from \cite{yakimov:leaves}. 
The $G_\Delta$-orbits in the wonderful
compactifications were also studied by Lusztig in \cite{lusz1} and
\cite{lusz2} using a somewhat different method.

\subsection{Some results on Weyl groups and generalized BD-triples}
\label{sec_weyl}

The results in this section, while different in presentation,
are closely related to some combinatorial results in \cite{bedard-comb},
which were used and extended in \cite{lusz1} and \cite{lusz2}.
In particular, the limit of a sequence studied in \cite{bedard-comb}
is closely related to the set $S(v, d)$ in 
Proposition \ref{prop_sequences}.

\begin{nota}
\label{nota_main-2}
{\em
Let $W$ be the Weyl group of $\Gamma$.
For $F \subset \Gamma$, let $W_F$ 
be the subgroup of $W$ generated by elements
in $F$. If $E, F \subset \Gamma$, let
${}^{E} \! W^{F}$ be the set of minimal
length representatives for  double cosets from $W_{E} \backslash W/W_{F}$,
and set $W^F= {}^{\emptyset}\!W^F$.
If $E_1, E_2 \subset F$, the set of
minimal length representatives in $W_F$ for the double cosets
from  $W_{E_1} \backslash W_F /W_{E_2}$ will be denoted by
${}^{E_1} \! (W_F)^{E_2}$. If $u \in
{}^{{E_1}}\! (W_F)^{{E_2}}$ and
$v  \in {}^{{E_{1}^{\prime}}}\! (W_{F^\prime})^{{E_{2}^{\prime}}}$,
we can regard both $u \in W_{F}$ and $v \in W_{F^\prime}$ as elements
in $W$, and by $uv$ we will mean their product in $W$.
}
\end{nota}
 
\begin{dfn}
\label{dfn_Svd}
{\em Let $(S, T, d)$ be a generalized BD-triple
in $\Gamma$. For $v \in W^T$, regarding $vd$ as a map $S \to \Delta$,
we define $S(v, d) \subset S$  to be the
largest subset in $S$ that is invariant under
$vd$. In other words,
\begin{equation}
\label{eq_Svd}
S(v,d)  = \{\alpha \in S: (vd)^n \alpha \in S, \forall
\, {\rm integer} \, \, n \geq 1\}.
\end{equation}
}
\end{dfn}

Parts 1) and 2) in the following Lemma \ref{lem_uw}
follow directly
from Proposition 2.7.5 of \cite{carter:lietype}
or Lemma 4.3 of \cite{yakimov:leaves}, and Part 3) is a special case of 
Lemma 5.3 in \cite{lu-milen1}.

\begin{lem}
\label{lem_uw}
1) If $w \in \swt$ and $u \in (W_S)^{S \cap w(T)}$, then $uw \in W^T$;

2) Every $v \in W^T$ has a unique decomposition $v = uw$,
where $w \in {}^S\!W^T$ and
$u \in (W_S)^{S \cap w(T)}$.
Moreover, $l(v) = l(u) + l(w)$;

3) For $w \in \swt$, set $T_w = S \cap w(T)$ and  $S_w = d^{-1}(T \cap w^{-1}(S))$
and regard $(S_w, T_w, wd)$ as a generalized BD-triple in $S$. 
Then for any $u \in (W_S)^{T_w}$, 
one has $S_w(u,wd) = S(uw, d)$, where
$S_w(u, wd)$ is the largest subset of $S_w$ that is
invariant under $uwd$.
\end{lem}

\begin{nota}
\label{nota_sequence}
{\em Fix a generalized BD-triple $(S, T, d)$. 
Let ${{\mathcal Q}}_{S, T, d}$ denote the set of all sequences ${\bf q} = 
\{{\bf q}_i\}_{i \geq 0}$ of quadruples 
\[
{\bf q}_i = (S_i,T_i,d_i,w_i), \hspace{.2in} i \geq 0,
\]
where, for each $i \geq 0$,

1) $(S_i,T_i,d_i)$ is a generalized BD-triple and $(S_0, T_0, d_0) =
(S, T, d)$;

2) $w_i \in {}^{S_i}\!(W_{S_{i-1}})^{T_i}$ (we set $S_{-1} = \Gamma$);

3) the triple $(S_{i+1}, T_{i+1}, d_{i+1})$ is obtained from 
${\bf q}_i$ as follows:
\[
T_{i+1} = S_i \cap w_i(T_i), \hspace{.2in} d_{i+1}=w_i d_i, \hspace{.2in}
S_{i+1} = d_{i+1}^{-1}(T_{i+1}).
\]
Note that $S_{i+1}, T_{i+1} \subset S_i$ for all $i$.
For ${\bf q}=\{{\bf q}_i\}_{i \geq 0} \in {{\mathcal Q}}_{S, T, d}$, let $i_0$
be the smallest integer such that $S_{i_0+1} = S_{i_0}$. Then it is easy to see 
that 
\[
{\bf q}_i = {\bf q}_{i_0+1} = (S_{i_0}, w_{i_0}(T_{i_0}), w_{i_0}d_{i_0}, 1)
\]
for all $i \geq i_0+1$, where $1$ is the identity element of $W$. 
Set 
\[
v_\infty({\bf q}) = w_{i_0} w_{i_0-1} \cdots w_0 
\hspace{.2in} {\rm and} \hspace{.2in}
S_\infty({\bf q}) = S_{i_0}.
\]
}
\end{nota}

Proposition \ref{prop_sequences}
is a direct consequence of Lemma \ref{lem_uw}. 
See also Proposition 2.5 of \cite{lusz1}.
 
\begin{prop}
\label{prop_sequences}
Let $(S, T, d)$ be a generalized BD-triple for $\Gamma$. Then for
any ${\bf q} \in {{\mathcal Q}}_{S, T, d}$, $v :=v_\infty({\bf q})
\in W^T$, and $S_\infty({\bf q}) = S(v, vd)$. Moreover, the map
${{\mathcal Q}}_{S, T, d} \rightarrow W^T: 
{\bf q} \mapsto v_\infty({\bf q})$
is bijective.
\end{prop}

\subsection{A double coset theorem}
\label{sec_a-double-coset-theorem}

For this section, $G$ will be a connected
complex reductive Lie group with Lie algebra $\fg$,
not necessarily of adjoint type. We use the same notation as
in Notation \ref{nota_main} and Notation \ref{nota_main-2}
for various subalgebras of
$\fg$ and subgroups of $G$.
We will define a class of subgroups $R$ of $G \times G$ that are slightly more
general than the groups $\RSTd$, and we will prove a theorem on 
$(G_\Delta, R)$-double cosets in $G \times G$ for such an $R$.

\begin{dfn}
\label{nota_main-3}
{\em
Let $(S, T, d)$ be a generalized BD-triple in $\Gamma$.
Let $C_S$ (resp. $C_T$) be a subgroup of the center $Z_S$ 
(resp. $Z_T$) of $M_S$ (resp. $M_T$), and let 
$\theta_d: M_S/C_S \to M_T /C_T$ be a group isomorphism
  that maps the one-dimensional
unipotent subgroup of $M_S/C_S$ defined by $\alpha$
to the corresponding subgroup of $M_T/C_T$ defined by 
$d\alpha$ for each $\alpha \in [S]$. 
By a $(S, T, d)$-admissible subgroup of $G \times G$ we mean
a subgroup $R = R(C_S, C_T, \theta_d)$ of $P_S \times P_T^-$ of the form
\begin{equation}
\label{eq_Q}
R(C_S, C_T, \theta_d) = \{(m, m^\prime) \in
M_S \times M_{T}: \, \theta_d(mC_S) = m^\prime C_T\} (N_S \times
N_{T}^{-}).
\end{equation}
}
\end{dfn}

\smallskip
Clearly $R(Z_S, Z_T, \gamma_d) = \RSTd$.
Let $R$ be any
$(S, T, d)$-admissible subgroup of $G \times G$.
Recall that the subset $S(v, d)$
of $S$  for $v\in W^T$ is defined in (\ref{eq_Svd}). If $\dot{v}$ is
a representative of $v$ in $G$, set
\[
R_{\dot{v}} =(M_{S(v,d)} \times M_{S(v,d)}) \cap
\left( ({\rm id} \times \Ad_{\dot{v}}) R\right),
\]
where $\Ad_{\dot{v}}: G \to G: g \mapsto \dot{v} g \dot{v}^{-1}$.
Let $R_{\dot{v}}$ act on $M_{S(v,d)}$ (from the right) by
\begin{equation}
\label{eq_Qdotv-on-MSvd}
m \cdot (m_1, m_1^\prime) = (m_1^\prime)^{-1} m m_1, \hs
m \in M_{S(v,d)}, \,
(m_1, m_1^\prime) \in R_{\dot{v}}.
\end{equation}
For $(g_1, g_2) \in G \times G$, let $[g_1, g_2]$ be the double coset $G_\Delta(g_1, g_2)R$ in $G \times G$.

\begin{thm}
\label{thm_double-cosets-Q}
Let $(S, T,d)$ be a generalized BD-triple, and
let $R = R(C_S, C_T, \theta_d)$ be an $(S, T, d)$-admissible
subgroup of $G \times G$
as given in (\ref{eq_Q}).
For $v \in W^T$, let
 $S(v,d) \subset S$ be given  in (\ref{eq_Svd}), and let
$\dot{v}$ be a fixed representative of $v$ in $G$.
Then

1) every $(G_\Delta, R)$-double coset in $G \times G$
is of the form $[m, \dot{v}]$ for some $v \in W^T$ and $m
\in M_{S(v, d)}$;

2) Two double cosets $[m_1, \dot{v}_1]$
and $[m_2, \dot{v}_2]$ in 1) coincide if and only if
$v_1 = v_2 = v$ and $m_1$ and $m_2$ are in the same $R_{\dot{v}}$-orbit
in $M_{S(v,d)}$ for the $R_{\dot{v}}$ action on $M_{S(v,d)}$ given
in (\ref{eq_Qdotv-on-MSvd}).
\end{thm}

We present the main induction step in the
proof of Theorem \ref{thm_double-cosets-Q} in a lemma.
Recall that each $w \in \swt$ gives rise to the 
generalized BD-triple $(S_w, T_w, wd)$ in $S$ as in
Lemma \ref{lem_uw}. Set
\[
N_{S_w}^{S} = N_{S_w} \cap M_S, \hs {\rm and} \hs
N_{T_w}^{S, -} = N_{T_w}^{-} \cap M_S.
\]
Fix a representative $\dw$ in $G$, and  define
\begin{equation}
\label{eq_QdwS}
R_{\dot{w}}^{S} =\left((M_{S_w} \times M_{T_w}) \cap   (({\rm id}
\times \Ad_{\dot{w}}) R)\right)
\left(N_{S_w}^{S} \times N_{T_w}^{S,-}\right) .
\end{equation}
Then $R_{\dot{w}}^{S}$ is an $(S_w, T_w, wd)$-admissible
subgroup of $M_S \times M_S$ defined by the subgroup
$C_S$ of $Z_{S_w}$, the subgroup $w(C_T)$ of $Z_{T_w}$ and the
group isomorphism
$\Ad_{\dot{w}} \theta_d:
M_{S_w}/C_S \to M_{T_w}/w(C_T)$.

\begin{lem}
\label{lem_double-cosets-1}
1) Every $(G_\Delta, R)$-double coset in $(G \times G)$ is of the form
$[m, m^\prime\dw]$ for a unique $w \in {}^S\!W^T$ and some
$m \in M_S $.

2) $[m_1, m_{1}^{\prime} \dw] = [m_2, m_{2}^{\prime} \dw]$, where $w \in {}^S\!W^T$
and $(m_1, m_{1}^{\prime}), (m_2, m_{2}^{\prime}) \in M_S \times M_S$,
if and only if
$(m_1, m_{1}^{\prime})$ and $(m_2, m_{2}^{\prime})$ are in the same
$((M_S)_\Delta, R_{\dot{w}}^{S})$-double coset in $M_S \times M_S$.
\end{lem}

\noindent
{\bf Proof.}
Consider the right action of $P_S \times P_{T}^{-}$ on $G_\Delta \backslash
(G \times G)$
by right translations. By the
Bruhat decomposition $G = \bigcup_{w \in {}^{S}\! W^T} P_S w P_{T}^{-}$,
 the set of $(P_S \times P_{T}^{-})$-orbits is parameterized by the set
$\{G_\Delta (e, \dw): w \in {}^S\! W^T\}$.
Let $w \in {}^S \! W^T$. The stabilizer subgroup of
$P_S \times P_{T}^{-}$ at $G_\Delta (e, \dw)$ is $P_S \cap (\dw P_{T}^{-}
\dw^{-1})$ considered as a subgroup of $\PSPT$ via the
embedding
\begin{equation}
\label{eq_PSwPT-embedding}
P_S \cap (\dw P_{T}^{-} \dw^{-1}) \lrw P_S \times P_{T}^{-}: \, \, p_S \Map
(p_S, \, \,
\dw^{-1}p_S \dw).
\end{equation}
Thus
the set of $R$-orbits in
$G_\Delta \backslash (G \times G)$ can be identified with
 the disjoint union over $w \in {}^S \! W^T$ of the spaces of
$R$-orbits in $P_S \cap (\dw P_{T}^{-} \dw^{-1}) \backslash
(P_S \times  P_{T}^{-})$.
Thus, for every $w \in {}^S \! W^T$, we have an injective map
\begin{equation}
\label{eq_orbitlemmaA}
(P_S \cap \dw P_{T}^{-} \dw^{-1}) \backslash
\PSPT/R \lrw G_\Delta \backslash G \times G/R
\end{equation}
given by
$(P_S \cap \dw P_{T}^{-} \dw^{-1})(p_S, p_{T}^{-})R
\to [p_S, \dw p_{T}^{-}]$.
We will complete the proof by identifying
\begin{equation}
\label{eq_orbitlemmaB}
(P_S \cap (\dw P_{T}^{-} \dw^{-1}) \backslash
\PSPT/R \; \cong \; (M_S)_\Delta \backslash M_S \times M_S /R_{\dot{w}}^S
 \end{equation}
through a series of steps.
Let $\pi_S: P_S \to M_S$ and $\pi_T: P_{T}^{-} \to M_T$ be the
projections with respect to the decompositions 
$P_S = M_S N_S$ and $P_T^- = M_T N_T^-$.  Then 
$\pi_S \times \pi_T: \PSPT \to M_S \times M_T$ gives an
identification
\begin{equation}
\label{eq_orbitlemmaC}
P_S \cap (\dw P_{T}^{-} \dw^{-1}) \backslash \PSPT/R \lrw
R_1
\backslash  M_S\times M_T /R_2,
\end{equation}
where $R_1 = (\pi_S \times \pi_T)(P_S \cap (\dw P_{T}^{-} \dw^{-1}))$ and $
R_2 =(M_S \times M_T) \cap R.$
Since the projection from $(M_S \times M_T) \cap R$ to $M_T$
is onto with kernel $(C_S \times \{e\})$, the map
\[
\phi_w: \, \, (M_S \times M_T)/ R_2 \lrw
(M_S \times M_S)/(M_S)_\Delta (C_S \times \{e\})
\]
that maps $(m_S, m_T)R_2$ to
$(m_{S}^{\prime}, m_S)((M_S)_\Delta (C_S \times \{e\}))$
is a well-defined bijection, where for $m_T \in M_T$,
$m_{S}^{\prime}$ is any element in $M_S$ such that
$(m_{S}^{\prime}, m_T) \in R_2.$
Thus $\phi_w$ induces an identification
\begin{equation}
\label{eq_orbitlemmaD}
\psi_w: \, \, R_1 \backslash M_S \times M_T /R_2
\lrw
R_3 \backslash M_S \times M_S /((M_S)_\Delta (C_S \times \{e\})),
\end{equation}
where
\[
R_3 \stackrel{{\rm def}}{=} \{(m_{S}^{\prime}, m_S) \in M_S \times M_S:
\, \exists
m_T \in M_T \, {\rm such} \, {\rm that} \,
(m_S, m_T) \in R_1, (m_{S}^{\prime}, m_T) \in R_2\}.
\]
By Theorem
2.8.7 of \cite{carter:lietype}, 
\begin{equation}
\label{eq_PScapPT}
P_S \cap (\dw P_{T}^{-} \dw^{-1}) =
(M_S \cap \Ad_{\dot{w}}(M_T))(M_S \cap \Ad_{\dot{w}}(N_{T}^{-}))
(N_S \cap \Ad_{\dot{w}}(M_T))(N_S \cap \Ad_{\dot{w}}(N_{T}^{-})).
\end{equation}
Note that $M_S \cap \Ad_{\dot{w}}(M_T)=M_{S \cap w(T)}$,
$M_S \cap \Ad_{\dot{w}}(N_{T}^{-}) = 
N^{S, -}_{S \cap w(T)} = M_S \cap N^{-}_{S \cap w(T)}$,
and $N_S \cap \Ad_{\dot{w}}(M_T) =
 N^{T}_{T \cap w^{-1}(S)}=M_T \cap N_{T \cap w^{-1}(S)}$.
Thus
\[
R_1 = \left\{(m, \Ad_{\dw^{-1}}(m)): m \in M_{S \cap w(T)}\right\} 
\left(N^{S, -}_{S \cap w(T)} \times N^{T}_{T \cap w^{-1}(S)}\right).
\]
Therefore $(m_{S}^{\prime}, m_S) \in R_3$ if and only if 
there exist $n \in N^{S, -}_{S \cap w(T)}$, $n_1 \in 
N^{T}_{T \cap w^{-1}(S)}$, and $m \in M_{S \cap w(T)}$ such that
$m_S = mn$ and $(m_{S}^{\prime}, \Ad_{\dw^{-1}}(m) n_1) \in R_2$.
It follows from the definition of $R$ that $(m_{S}^{\prime}, m_S) \in  R_3$
if and only if there exist $m^\prime \in M_{S_w}, m \in M_{T_w}$,
$n \in N^{S, -}_{T_w}$, and $n^\prime \in N^{S}_{S_w} = M_S \cap N_{S_w}$
such that $m_S = mn, m_{S}^{\prime} = m^\prime n^\prime$ and $(m^\prime,
\Ad_{\dot{w}}^{-1}(m)) \in R.$
Thus $R_3 = R_{\dot{w}}^{S}$.
Since 
$C_S \times \{e\} \subset R_{\dot{w}}^{S}$, the (right)
action of $C_S \times \{e\}$ on $R_{\dot{w}}^{S} \backslash
(M_S \times M_S)$ is trivial. Thus we have
\beqa
\label{eq_orbitlemmaE}
R_3\backslash M_S \times M_S /((M_S)_\Delta (C_S \times \{e\}))
 & \cong &
R_{\dot{w}}^{S} \backslash M_S \times M_S /((M_S)_\Delta (C_S \times \{e\}))
\\
&  \cong &
R_{\dot{w}}^{S} \backslash M_S \times M_S /(M_S)_\Delta \\
& \cong &(M_S)_\Delta \backslash M_S \times M_S /R_{\dot{w}}^S,
\eeqa
where the last identification is induced by the inverse map
of $M_S \times M_S$. 

Combining the above identification with the
identifications in (\ref{eq_orbitlemmaC})-(\ref{eq_orbitlemmaD}) and 
the inclusion of (\ref{eq_orbitlemmaA}), we
get a well-defined injective map
$(M_S)_\Delta \backslash M_S \times M_S /R_{\dot{w}}^S \to
G_\Delta \backslash G \times G /R$
given by 
\[
(M_S)_\Delta (m,m^\prime) {R}_{\dot{w}}^S \to
[((m^\prime)^{-1},\dot{w} \theta_d(m^{-1}))] =
[(m^{\prime})^{-1}m, \dot{w}]=[m,m^{\prime}\dot{w}].
\]
 This finishes the proof of
 Lemma \ref{lem_double-cosets-1}.
\qed

% Let $\Gamma=S_{i-1}$, let $S=S_i$, $d=d_i$, and $T=T_i$.
% Let $w\in {}^{S_i}W_{S_{i-1}}^{T_i}$. Then 
% $R_i=R_{\dot{w}}^{S}$.

\smallskip
\noindent
{\bf Proof of Theorem \ref{thm_double-cosets-Q}.}
By Lemma \ref{lem_double-cosets-1}, each $(G_\Delta,R)$
double coset in $G\times G$ determines a unique
$w\in {}^{S}\!W^{T}$ and a unique
double coset $[m,m^\prime]_1 \in 
(M_{S})_\Delta \backslash M_{S} \times M_{S} / R_{\dot{w}}^{S}$.
Let $(S_0, T_0, d_0, w_0) = (S, T, d, w)$.
By successively applying Lemma \ref{lem_double-cosets-1}
to a sequence of smaller subgroups,
we obtain a sequence ${\bf q}$ of quadruples  ${\bf q}_i = (S_i,T_i,d_i,w_i)$
as in Notation \ref{nota_sequence}, as well as a double coset
in $(M_{S_i})_\Delta \backslash M_{S_i} \times M_{S_i} / R_i,$
where $R_i$ is the subgroup of $M_{S_i} \times M_{S_i}$
defined analogously to $R_{\dot{w}}^{S}$.

As in  Notation \ref{nota_sequence}, let $i_0$ be the smallest
integer such that $S_{i_0+1} = S_{i_0}$ and let
$v = v_\infty({\bf q})=w_{i_0} w_{i_0-1} \cdots w_0$.
Then each $(G_\Delta, R)$-double coset in $G \times G$ is of the form 
$[m,m^{\prime}\dot{v}]$
 for $m\in M_{S_{i_0+1}}$. 
By
Proposition \ref{prop_sequences}, $v \in W^T$, and $S_{i_0+1} = S(v, d)$.
Moreover,  $R_{i_0+1} = R_{\dot{v}}$ by definition. Thus 
double cosets
in 
$(M_{S_{i_0+1}})_\Delta \backslash 
M_{S_{i_0+1}} \times M_{S_{i_0+1}} / R_{i_0+1}$
coincide with double cosets in
$(M_{S(v,d)})_\Delta \backslash
M_{S(v,d)} \times M_{S(v,d)} / R_{\dot{v}}$.
It is easy to see that the map
\[
(M_{S(v,d)})_\Delta \backslash
M_{S(v,d)} \times M_{S(v,d)} / R_{\dot{v}} \lrw M_{S(v,d)}/
 R_{\dot{v}}: \, \,
[m, m^\prime] \Map [{m^{\prime}}^{-1}m]
\]
is a bijection. This proves Theorem \ref{thm_double-cosets-Q}.
\qed

\subsection{$G_\Delta$-orbits in $\lag$}
\label{sec_G-orbits}

\begin{nota}
\label{nota_G-orbits-L}
{\em
For a generalized BD-triple $(S, T, d)$,
$V \in
\LinLagr(\fz_S \oplus \fz_T)$, $m \in {M_{S(v,d)}}$,
$v \in W^T$, and $\dot{v} \in G$ a fixed representative of $v$ in $G$,
set
\begin{equation}
\label{eq_flSTdVvm}
\fl_{S, T, d, V, \dot{v}, m}  = \Ad_{(m, \dot{v})}
\fl_{S, T, d, V},
\end{equation}
where $\fl_{S, T, d, V}$ is given in (\ref{eq_flSTdV}).
Define
\begin{equation}
\label{eq_Rdotv-2}
R_{\dot{v}} = \{(m_1, m_1^\prime) \in M_{S(v,d)}
\times M_{S(v,d)}: \, \gamma_d (\chi_S(m_1)) = \chi_T(\dot{v}^{-1}
m_1^\prime \dot{v})\},
\end{equation}
and let $R_{\dot{v}}$ act on $M_{S(v,d)}$ (from the right) by
\begin{equation}
\label{eq_Rdotv-on-MSvd}
m \cdot (m_1, m_1^\prime) = (m_1^\prime)^{-1} m m_1, \hs
m \in M_{S(v,d)}, \,
(m_1, m_1^\prime) \in R_{\dot{v}}.
\end{equation}
}
\end{nota}

As an immediate corollary of Theorem \ref{thm_double-cosets-Q}, we have
\begin{cor}
\label{cor_G-orbits-L}
Every $G_\Delta$-orbit
in $\lag$ passes through an $\fl_{S, T, d, V, \dot{v}, m}$ for a
unique generalized BD-triple $(S, T, d)$,
a unique $V \in
\LinLagr(\fz_S \oplus \fz_T)$, a unique $v \in W^T$, and
some $m \in {M_{S(v,d)}}$; Two such
Lagrangian subalgebras $\fl_{S, T, d, V, v, m_1}$
and $\fl_{S, T, d, V, v, m_2}$
are in the same $G_\Delta$-orbit if and only if
$m_1$ and $m_2$ are
in the same $R_{\dot{v}}$-orbit in $M_{S(v,d)}$.
\end{cor}

\subsection{Normalizer subalgebras of $\fg_\Delta$ at $\fl \in \lag$}
\label{sec_normalizers}

For $\fl =\fl_{S, T, d, V, \dot{v}, m}$ as
in Corollary
\ref{cor_G-orbits-L}, we now compute its normalizer subalgebra
$\fn(\fl)$ in $ \fg \cong
\fg_\Delta = \{(x, x): x \in \fg\}$.
Introduce
\[
\phi :=\Ad_{\dot{v}} \gamma_d \chi_S \Ad_{m}^{-1}: \, \, \, \fp_S \lrw \fg.
\]
Consider the standard parabolic subalgebra $\fp_{S(v,d)}$ and its
decomposition $\fp_{S(v,d)} = \fz_{S(v,d)} + \fg_{S(v,d)} + \fn_{S(v,d)}$
(see Notation \ref{nota_main}).

\begin{lem}
\label{lem_phi-three-pieces}
The map $\phi = \Ad_{\dot{v}} \gamma_d \chi_S \Ad_{m}^{-1}$ leaves each of
$\fz_{S(v,d)}$, $\fg_{S(v,d)}$, and  $\fn_{S(v,d)}$ invariant. Moreover,
$\phi: \fn_{S(v,d)} \to \fn_{S(v,d)}$ is nilpotent.
\end{lem}

\noindent
{\bf Proof.} Let $x \in \fz_{S(v,d)}$. Then $\phi(x) = \Advd \chi_S(x) \in
\fh$. For $\alpha \in S(v,d)$, since $(vd)^{-1}\alpha \in S(v,d)$,
\[
\alpha(\phi(x)) = ((vd)^{-1}\alpha)(\chi_S(x)) = ((vd)^{-1}\alpha)(x) =0.
\]
Thus $\phi(x) \in \fzsvd$, so $\fzsvd$ is $\phi$-invariant.
Since both $\Advd$ and $\Ad_{m}^{-1}$ leave $\fgsvd$ invariant, we see that
$\phi|_{{\mathfrak g}_{S(v,d)}}= \Advd \Ad_{m}^{-1}$ leaves $\fgsvd$ invariant.

To show that $\fnsvd$ is $\phi$-invariant and that
$\phi: \fn_{S(v,d)} \to \fn_{S(v,d)}$ is nilpotent, set
$\Sigma_{0}^{+} = \Sigma^{+} -[S]$, and for $j \geq 1$, set
\begin{equation}
\label{eq_Sigmaj}
\Sigma_{j}^{+} =  \{\alpha \in \Sigma^+: \, \alpha \in [S], 
vd\alpha \in [S], \cdots,
(vd)^{j-1}(\alpha) \in [S], (vd)^j(\alpha) \notin [S]\}.
\end{equation}
Then $\Sigma^+ -[S(v,d)] = \cup_{j \geq 0} \Sigma_{j}^{+}$ and
$vd(\Sigma_j^+) \subset  \Sigma_{j-1}^{+}$
for $j \geq 1$.
For $j \geq 0$, set
$\fn_j = \oplus_{\alpha \in \Sigma_{j}^{+}} \fg_\alpha.$
Then $\fn_0 = \fn_S$, and
$\fnsvd = \sum_{j \geq 0} \fn_j$
is a finite direct sum. 
It is easy to prove by induction on $j$ that
\begin{equation}
\label{eq_alpha-beta}
\alpha \in \Sigma_{j}^{+}, \beta \in [S(v,d)], \alpha + \beta \in \Sigma \, \,
\Longrightarrow \, \alpha + \beta \in \Sigma_{j}^{+}, \hspace{.2in} \forall j \geq 0.
\end{equation} 
It follows 
that $[\fmsvd, \fn_j]
\subset \fn_j$ for each $j \geq 0$. Thus
$\Ad_{m} \fn_j = \fn_j, \forall j \geq 0.$
By setting $\fn_{-1} = 0$, we then have 
$ \Ad_{\dot{v}} \gamma_d \chi_S(\fn_j) \subset \fn_{j-1},  \forall j \geq 0$. Thus
$\phi(\fn_j) \subset \fn_{j-1}$ for all $j \geq 0$, and 
$\phi: \fn_{S(v,d)} \to \fn_{S(v,d)}$ is nilpotent.
\qed

\begin{rem}
\label{rem_minuscase}
{\em
Similar arguments imply the same statements
 for $\fn^{-}_{S(v, d)}$. In particular,
for $j\ge 0$, 
set $\fn_j^- = \oplus_{\alpha \in \Sigma_{j}^{+}} \fg_{-\alpha}$
and set $\fn_{-1}^- = 0$.
Then $\phi(\fn_j^-)\subset \fn_{j-1}^-$ for $j\ge 0$,
$\fn_S^- = \fn_0^-$, and $\fn_{S(v,d)}^- = \sum_{j\ge 0} \fn_j^-$
is a direct sum.
}
\end{rem}

Since $\phi: \fnsvd \to \fnsvd$ is nilpotent, we can define
\[
\psi:= (1-\phi)^{-1} = 1 + \phi + \phi^2 + \phi^3 + \cdots: \, \,
\fnsvd \lrw \fnsvd.
\]
Let $\Sigma_v^+ = \{\alpha \in \Sigma^{+}: v^{-1} \alpha \in \Sigma^-\}$.
Since $v ([T] \cap \Sigma^+)
\subset \Sigma^+$,  we have 
$\Sigma_v^+ \subset \Sigma^+ -[S(v,d)]$. Let
\[
\fn_v = \oplus_{\alpha \in \Sigma_v^+} \fg_\alpha = \fn \cap \Adv (\fn^-).
\]
Then $\fn_v \subset \fnsvd$.

\begin{thm}
\label{thm_normalizers}
The normalizer subalgebra $\fn(\fl)$ in $\fg_\Delta \cong \fg$
of $\fl = \fl_{S, T, d, V, \dot{v}, m}$
 in (\ref{eq_flSTdVvm}) is
\[
\fn(\fl) = \fz_{S(v,d)}^{\prime} + \fg_{S(v, d)}^{\phi} + \psi(\fn_v),
\]
where $\fg_{S(v, d)}^{\phi}$ is
the fixed point set of $\phi|_{{\mathfrak g}_{S(v,d)}} = \Advd \Ad_{m}^{-1}$
 in $\fgsvd$, and
\[
\fz_{S(v,d)}^{\prime}=\{z \in \fzsvd: z - \phi(z) \in \Adv \fz_T\}
 = \{z \in \fzsvd: \gamma_d \chi_S(z)=\chi_T (\Ad_{\dot{v}}^{-1} z)\}.
\]
\end{thm}

\noindent
{\bf Proof.} The normalizer subgroup $\RSTd$ of $\fl_{S, T, d, V}$ in 
$G \times G$ has Lie algebra
\[
\fr_{S, T, d}= (\fz_S \oplus \fz_T) + (\fn_S \oplus \fn_{T}^{-}) +
\{(x, \gamma_d(x)): x \in \fg_S\} = 
\{(x,y)\in \fp_S \oplus \fp_T^- : \gamma_d \chi_S(x) =\chi_T(y) \}.
\]
 Since $\fl = \Ad_{(m,\dv)}\fl_{S, T, d, V}$, it follows that
 $\fn(\fl) = \{x \in \fg: \, (\Ad_{m}^{-1} x, \, \Ad_{\dot{v}}^{-1} x) 
\in \fr_{S, T, d}\}$.
 Thus  $x \in \fn(\fl)$ if and only if $x \in \fp_S \cap
\Adv \fp_{T}^{-}$ and
$\gamma_d \chi_S(\Ad_{m}^{-1}(x)) = \chi_T(\Ad_{\dot{v}}^{-1} x)$, which is
equivalent to 
\begin{equation}
\label{eq_for-x}
x - \Advd  \chi_S(\Ad_{m}^{-1}(x)) \in \Adv (\fz_T+\fn_{T}^{-}).
\end{equation}
Let $\chi_S$ also denote the projection $ \fg \to \fg_S$  with respect to
the decomposition
$\fg = \fn_{S}^{-} + \fz_S + \fg_S + \fn_S,$ so
$\phi: x \mapsto \Advd \chi_S(\Ad_{m}^{-1}x)$ is defined for all $x \in \fg$.
Let $\fc$ be the  set of all $x \in \fg$ satisfying (\ref{eq_for-x}). 
Then $\fn(\fl) =\fc \cap 
\left(\fp_S \cap \Adv \fp_{T}^{-}\right)$.

Consider the decomposition
$\fg = \fnsvdp + \fmsvd + \fnsvd$. 
Since $\alpha \notin v([T])$ implies that
$\alpha \notin [S(v, d)]$, we have
\[
\Adv \fn_{T}^{-} = \left(\Adv \fn_{T}^{-} \right) \cap \fn^- +
 \left(\Adv \fn_{T}^{-}\right)\cap \fn 
\; \subset \;  \fnsvdp + \fnsvd.
\]
Moreover, it is easy to see that $\left(\Adv \fn_{T}^{-}\right)\cap \fn = 
\left(\Adv \fn_{T}^{-}\right)\cap \fn_{S(v,d)} =\fn_v$, so  
\begin{equation}
\label{eq_AdvnT}
\Adv \fn_{T}^{-} = \left(\Adv \fn_{T}^{-} \right) \cap \fn^{-}_{S(v,d)} + \fn_v.
\end{equation}
 Now let $x \in \fg$ and write $x = x_- + x_0 + x_+$, where
$x_- \in \fnsvdp, x_0 \in \fmsvd,$ and $x_+ \in \fnsvd$.
It follows from Lemma \ref{lem_phi-three-pieces} and
(\ref{eq_AdvnT}) that $x \in \fc$, i.e., 
$x$ satisfies (\ref{eq_for-x}), if and only if
 \begin{equation}
 \label{eq_x-fc}
 \begin{cases} x_0 - \phi(x_0) \in \Adv \fz_T \\
 x_+ - \phi(x_+) \in \fn_v\\
 x_- - \phi(x_-) \in \left(\Adv \fn_{T}^{-} \right) \cap \fn^{-}_{S(v, d)}.
 \end{cases}
\end{equation}
Write $x_0 = z_0 + y_0$, where $z_0 \in \fzsvd$
and $y_0 \in \fgsvd$. Since $\Adv \fz_T \subset \fzsvd$ and since
both $\fzsvd$ and $\fgsvd$ 
are $\phi$-invariant, $x_0 - \phi(x_0) \in \Adv \fz_T$ if and only if
$z_0 -\phi(z_0) \in \Adv \fz_T$ and $y_0 - \phi(y_0) = 0$, which 
is the same as $x_0 \in \fz_{S(v,d)}^{\prime} + \fg_{S(v,d)}^{\phi}$.
Recall that $\psi=(1-\phi)^{-1}$ on $\fnsvd$.  Thus,
$x_+ - \phi(x_+) \in \fn_v$ if and only if $x_+ \in \psi(\fn_v)$.
 Since $\fn_v \subset \fn_{S(v,d)}$,
\[
\psi(\fn_v)
\subset \fp_S\cap ( \fn_v + \phi(\fnsvd)) 
\subset \fp_S \cap (\Adv(\fn^-  +  \fm_T)) \subset  
\fp_S \cap \Adv \fp_{T}^{-}.
\]
Note that $\fz_{S(v,d)}^{\prime} + \fg_{S(v, d)}^{\phi} \subset \fp_S 
\cap \Adv \fp_{T}^{-}$.
Thus $\fc \cap (\fp_S \cap \Adv \fp_{T}^{-}) = \fz_{S(v,d)}^{\prime} + 
\fg_{S(v, d)}^{\phi} +\psi (\fn_v) +\fc^\prime,$
where $\fc^\prime$ consists of all 
\[
x_- \in \fn_{S(v, d)}^- 
\cap\fp_S \cap \Adv \fp_{T}^{-} \subset \fn_{S(v, d)}^- 
\cap\fm_S
\]
satisfying the third condition in (\ref{eq_x-fc}). It suffices to show that
$\fc^\prime = 0$.

We regard the direct sum decomposition  
$\fn_{S(v,d)}^- = \sum_{j\ge 0} \fn_j^-$ 
from Remark \ref{rem_minuscase} as a grading of 
$\fn_{S(v,d)}^-$. Let 
$U = \fn_{S(v, d)}^- \cap \fm_S = \sum_{j > 0} \fn_j^-$
and let $Y= \Adv \fn_{T}^{-} \cap \fn^{-}_{S(v, d)}$.
Clearly, $U$ and $Y$ are graded subspaces of $\fn^{-}_{S(v,d)}$,
since they are sums of root spaces.
Since $U\subset \fg_S$, $\phi$ is injective on $U$. Moreover,
the image of $\phi$ is in $\Adv(\fg_{T})$ so $Y$ has zero
intersection with the image of $\phi$. The fact that $\fc^\prime = 0$ now
follows from the following simple linear algebra fact in Lemma \ref{lem_linalgfact}.
\qed

\begin{lem}
\label{lem_linalgfact}
Let $V=\oplus V_i$ be a graded vector space with graded subspaces
$U$ and $Y$. Let $\phi$ be an endomorphism of $V$ such that
1) $\phi(V_i)\subset V_{i-1}$ for all $i$,
2) $Y \cap  {\rm Im}(\phi) = 0$, and 
3) $U \cap {\rm Ker}(\phi) = 0$.
Then $\{ v \in U : v - \phi(v) \in Y \} = 0$.
\end{lem}

\begin{rem}
\label{rem_why-subalgebra}
{\em   
Theorem \ref{thm_normalizers} implies that
\[
\fn(\fl) \subset \fp_{S(v,d)} \cap \Adv \fp_{T}^{-}
=\fp_{S(v,d)} \cap \Adv \fp_{T}^{-} = \fmsvd + \fn_v +
\fnsvd \cap \Adv \fm_T
\]
and that $\fn(\fl) = \fn(\fl) \cap \fmsvd + 
\fn(\fl) \cap (\fn_v +\fnsvd \cap \Adv \fm_T)$, where
\[
(\fn(\fl) \cap \fmsvd)_\Delta = 
(\fz_{S(v,d)}^{\prime} + \fg_{S(v, d)}^{\phi})_\Delta 
= (\fmsvd)_\Delta \cap \Ad_{(m, \dot{v})} \fr_{S, T, d},
\]
and $\fn(\fl) \cap (\fn_v +\fnsvd \cap \Adv \fm_T)$ is the graph of the map
\[
\psi - 1 = \phi \psi = \phi + \phi^2 + \cdots: 
\, \fn_v  \to \fnsvd \cap \Adv \fm_T.
\]
In \cite{lu-milen1}, the map $\psi-1$ is shown to be related to some
set-theoretical solutions
to the Quantum Yang-Baxter Equation.
} \end{rem}

 \subsection{Intersections of $\fg_\Delta$ with  $\fl \in
{{\mathcal L}}$}
\label{sec_intersections}

By Corollary
\ref{cor_G-orbits-L}, to compute $\fg_\Delta \cap \fl$ for any $\fl \in
{{\mathcal L}}$, we may assume that 
$\fl = \fl_{S, T, d, V, \dot{v}, m}$
as given in (\ref{eq_flSTdVvm}).

\begin{prop}
\label{prop_intersections}
For the Lagrangian subalgebra $\fl_{S, T, d, V, \dot{v}, m}$
as given in (\ref{eq_flSTdVvm}), let the notation be as 
in Theorem \ref{thm_normalizers}.
Then
\[
\fg_\Delta \cap  \fl_{S, T, d, V, \dot{v}, m} = \Ad_{(m,\dot{v})}V^\prime
+ \left(\fg_{S(v,d)}^{\phi} + \psi(\fn_v)\right)_\Delta,
\]
where $V^\prime = \{(z,  { {v}}^{-1}z): z \in \fz_{S(v,d)}^{\prime}\}
\cap \left(V + \{(x, \gamma_d(x)): x \in \fh_S\}\right)$.
\end{prop}

\noindent
{\bf Proof.} Set $\fl = \fl_{S, T, d, V, \dot{v}, m}$.
By Theorem \ref{thm_normalizers}, 
\[
\fg_\Delta \cap \fl \subset \fn(\fl) 
=\left(\fz_{S(v,d)}^{\prime} + \fg_{S(v,d)}^{\phi} + \psi(\fn_v)\right)_\Delta.
\]
Since $\left(\fg_{S(v,d)}^{\phi} + \psi(\fn_v)\right)_\Delta 
\subset \fl$, we see that
\[
\fg_\Delta \cap \fl = \left((\fz_{S(v,d)}^{\prime})_\Delta \cap \fl\right) +
\left(\fg_{S(v,d)}^{\phi} + \psi(\fn_v)\right)_\Delta,
\]
and
\[
\left(\fz_{S(v,d)}^{\prime}\right)_\Delta \cap \fl =
  \left(\fz_{S(v,d)}^{\prime}\right)_\Delta \cap \fl \cap (\fh \oplus \fh)
= \Ad_{(m,\dot{v})}V^\prime.
\]
\qed

Recall that a Belavin-Drinfeld triple \cite{b-d:factorizable}
for $\fg$ is a generalized Belavin-Drinfeld triple $(S, T, d)$ 
with the nilpotency condition:
for every $\alpha \in S$, there
exists an integer $n \geq 1$ such that $\alpha, d\alpha, \ldots, d^{n-1}
\alpha \in S$ but $d^n \alpha \notin S$. The nilpotency condition is
equivalent to  $S(1, d) = \emptyset$, 
where $1$ is the
identity element in the Weyl group $W$.

\begin{dfn}
\label{dfn_BD-system}
{\em A {\it Belavin-Drinfeld system} is 
a quadruple $(S, T, d, V)$, where $(S, T, d)$ is 
a Belavin-Drinfeld triple, and 
$V$ is a Lagrangian subspace of $\fz_S \oplus \fz_T$ such that
\[
\fh_\Delta \cap \left(V + \{(x, \gamma_d(x)): \, x \in \fh_S\}\right) = 0.
\]
}
\end{dfn}

We now derive a theorem of Belavin and Drinfeld \cite{b-d:factorizable}
from Proposition \ref{prop_intersections}.

\begin{cor}
\label{cor_BD} [Belavin-Drinfeld]
A Lagrangian subalgebra $\fl$ of $\fg \oplus \fg$ has trivial 
intersection with
$\fg_\Delta$ if and only if $\fl$ is $G_\Delta$-conjugate to a Lagrangian
subalgebra of the form $\fl_{S, T, d, V}$, where $(S, T, d, V)$ 
is a Belavin-Drinfeld
system.
\end{cor}

\noindent
{\bf Proof.} It is clear from Proposition \ref{prop_intersections}
that $\fg_\Delta \cap \fl_{S, T, d, V} = 0$ if $(S, T, d, V)$ is a
Belavin-Drinfeld system. Suppose that 
$\fg_\Delta \cap \fl_{S, T, d, V, \dot{v}, m} = 0$, where 
$\fl_{S, T, d, V, \dot{v}, m}$
is as in (\ref{eq_flSTdVvm}).
Since $\dim \psi(\fn_v) = l(v)$,
the length of $v$, and since every automorphism of a
semi-simple Lie algebra has 
fixed point set of dimension at least one \cite{wi:fixedset},
  $v=1$
and $S(1, d) = \emptyset$.
In this case, $V^\prime$ as in Proposition  \ref{prop_intersections} is
given by
\[
V^\prime = \fh_\Delta \cap \left(V + \{(x, \gamma_d(x)): \, x \in \fh_S\}\right),
\]
so $\fh_\Delta \cap \left(V + \{(x, \gamma_d(x)): \, x \in \fh_S\}\right) =0$, and
we have $\fl_{S, T, d, V, \dot{v}, m} = \Ad_{(m, \dot{v})} \fl_{S, T, d, V}$
for some $m \in H$ and $\dot{v} \in H$. Note that in this case
\[
R_{\dot{v}} = \{(h_1, h_2) \in H \times H: \gamma_d(\chi_S(h_1)) = 
\chi_T(h_2)\}
\]
and
$R_{\dot{v}}$ acts on $H$ from the right by
$h \cdot(h_1, h_2) = h h_1 h_{2}^{-1}$, where $h \in H$ and 
$(h_1, h_2) \in R_{\dot{v}}$. Consider
${\bf m}:  R_{\dot{v}} \to H:  (h_1, h_2) \mapsto h_1 h_{2}^{-1}.$
The assumption $V^\prime = 0$ implies that the dimension of the
kernel of the differential of ${\bf m}$ is less than or equal to
$\dim(\fz_T)$. It follows that the differential of ${\bf m}$
is onto, so 
${\bf m}$ is onto. By Corollary \ref{cor_G-orbits-L},
$\fl_{S, T, d, V, \dot{v}, m}$ is in the   $G_\Delta$-orbit of
$\fl_{S, T, d, V}$.
 \qed

\subsection{Examples of smooth $G_\Delta$-orbit closures in $\lag$}
\label{G-Delta-closures}

The closure of a $G_\Delta$-orbit in $\lag$ is  not
necessarily smooth. We now look at two cases for which such a
closure is smooth.

\begin{prop}
\label{prop_G-closure-lBD}
If $\fl \in \lag$ is such that
$\fg_\Delta \cap \fl = 0$, then $\overline{G_\Delta \cdot \fl} = 
\overline{(G \times G)\cdot \fl}$ is smooth.
\end{prop}

\noindent
{\bf Proof.} We only need to show that $\dim(G_\Delta \cdot \fl)=
\dim((G \times G)\cdot \fl)$.
By Corollary \ref{cor_BD}, we may assume that
$\fl = \fl_{S,T, d, V}$, where $(S, T, d, V)$ is a Belavin-Drinfeld
system, and so  
\[
\fg_\Delta \cap \fr_{S, T, d} = \fh_{\Delta} \cap
((\fz_S \oplus \fz_T) + V_S)),
\]
where $V_S = \{(x, \gamma_d(x)): x \in \fh_S\}$. For a subspace
$A$ of $\fh \oplus \fh$, let
\[
A^\perp = \{(x, y) \in \fh \oplus \fh: \,
\la (x, y), (x_1, y_1) \ra = 0 \, \forall
(x_1, y_1) \in A\}.
\]
 Then
$\left(\fh_{\Delta} \cap
((\fz_S \oplus \fz_T) + V_S))\right)^\perp = \fh_\Delta + V_S.$
Since $\fh_\Delta \cap V_S = 0$, we have 
\[
\dim(\fh_{\Delta} \cap
((\fz_S \oplus \fz_T) + V_S))) = 2\dim \fh - \dim \fh - \dim \fh_S =\dim \fz_S. 
\]
Thus
$\dim (G_\Delta \cdot \fl) =\dim \fg - \dim \fz_S= \dim ((G \times G)\cdot \fl)$
by Proposition \ref{prop_GG-orbit-types}.
\qed

We now show that the De Concini-Procesi compactification of a complex
symmetric
space of $G$ can be embedded into $\lag$ as the closure
of a $G_\Delta$-orbit in $\lag$.

Let $\sigma:\fg\to\fg$ be an involution with lift $\sigma$ to $G$, and
let
$\fg^\sigma$ and $G^\sigma$ be the fixed subalgebra and subgroup of $\sigma.$
Let again $\fl_{\sigma} \in \lag$ be the graph of $\sigma$. The orbit
 $G_\Delta \cdot \fl_\sigma$ may be identified with the complex
symmetric space $G/G^\sigma.$ 
We will
show that the closure $\overline{G_\Delta \cdot \fl_\sigma}$
is isomorphic to 
the De Concini-Procesi compactification of $G/G^\sigma,$ which
as defined as follows. Let $\dim(\fg^\sigma)=m$, so
$\fg^\sigma \in {\rm Gr} (m,\fg).$ Then $G\cdot \fg^\sigma\cong
G/G^\sigma,$ and $X_\sigma := \overline{G\cdot \fg^\sigma}$, the closure of
$G\cdot \fg^\sigma$ in $ {\rm Gr} (m,\fg)$, is
the De Concini--Procesi compactification. It is smooth
with finitely many $G$-orbits \cite{dp:compactification}.

We recall some basic results about involutions. 
Choose a
$\sigma$-stable maximal split Cartan subalgebra $\fhsplit$ of $\fg$,
i.e., a $\sigma$-stable Cartan subalgebra $\fhsplit$ such that
${\fh}_s^{-\sigma}$ has maximal dimension.
There is an induced action of $\sigma$ on
the roots of $\fhsplit$ in $\fg$, and there is a positive root 
system $\Phiplushs$
 for $\fhsplit$
with the property that if $\alpha\in\Phiplushs$, then either 
$\sigma(\alpha)=\alpha$ and $\sigma|_{\fg_\alpha} = \id$,
 or $\sigma(\alpha)\not\in\Phiplushs$.
A weight $\lambda \in {\fh}_s^*$ is called a regular special dominant
weight
if $\lambda$ is nonnegative on roots in $\Phiplushs$, 
$\sigma(\lambda)=-\lambda$, and $\lambda(H_\alpha)=0$ for $\alpha$
simple implies that $\sigma(\alpha)=\alpha$.
If $\lambda$ and $\mu$ are weights, we say $\lambda \ge \mu$
if $\lambda - \mu = \sum_{\alpha \in \Phiplushs, n_\alpha \ge 0}
 n_\alpha \alpha$. For a weight $\mu$, let $\overline{\mu} =
\frac{1}{2}(\mu - \sigma(\mu))$.

\begin{lem}
\label{lem_dplemma}[De Concini-Procesi,  \cite{dp:compactification},
Lemmas 4.1 and 6.1]
Let $V$ be a representation of $G$, and suppose there exists
a vector $v\in V$ such that $G^\sigma$ is the stabilizer of
the line through $v$. Suppose that when we decompose 
$v$ into a sum of weight vectors for $\fhsplit$,
$v=v_\lambda + \sum v_i$ where $v_\lambda$ has regular
special dominant weight $\lambda$ and each $v_i$ has weight
$\mu_i$ where $\lambda \ge \overline{\mu_i}$. Let
$[v]$ be class of $v$ in $\Proj(V)$ 
and let $X^\prime$ be the closure of $G\cdot [v]$ in
$\Proj(V)$. Then $X^\prime \cong
X_\sigma$.
\end{lem}

\begin{prop}
\label{prop_orbitequalsdp}
There is a $G$-equivariant isomorphism $\overline{\GDelta\cdot \flsigma}
\cong X_\sigma.$
 \end{prop}

\noindent{\bf Proof.}
To apply Lemma \ref{lem_dplemma}, let $n=\dim(\fg)$ and consider the
diagonal action
of $G$ on $V=\wedge^n(\fg\oplus\fg)$ and the vector
$v_\sigma =\wedge^n(\flsigma).$
In order to represent $v_\sigma$ as a sum of weight vectors
in $ \wedge^n(\fg\oplus
\fg),$ we choose a basis.
Let $U_1,\dots, U_l$ be a basis of $\fhsplit$.
Let $\beta_1, \dots, \beta_s$ be the roots of $\Phiplushs$
 such that $\sigma(\beta_i)=\beta_i$, and let $\alpha_1, \dots, \alpha_t$ be the other roots in
$\Phiplushs$. For each root $\alpha,$ choose a root vector $X_\alpha.$
Then
\[
\{ (U_i,\sigma(U_i))| i=1,\dots, l \} \cup
\{ (X_{\pm\beta_i}, X_{\pm\beta_i}) | i=1, \cdots, s \} \cup
\{ (X_{\pm\alpha_j}, \sigma(X_{\pm\alpha_j})) | i=1, \cdots t \}
\]
is clearly a basis of $\flsigma$.
Now $v_\sigma$ is the wedge of the vectors $(Y_i,\sigma(Y_i))$ as $Y_i$
runs through the above basis, and $v_\sigma$ contains the summand
\[
u:\bigwedge_{i=1,\dots, l} (U_i,\sigma(U_i))
\bigwedge_{i=1, \dots, s} (X_{\beta_i},0)\wedge (X_{-\beta_i},0)
\bigwedge_{j=1, \dots, t} (X_{\alpha_i},0)\wedge (0,\sigma(X_{-\alpha_i})).
\]
It is easy to see that
$u$ is a weight vector for the diagonal Cartan subalgebra with weight
$\nu:=\sum_{i=1, \dots, t} \alpha_i - \sigma(\alpha_i),$ and
 $\nu = 2\sum_{i=1,\dots, t} \alpha_i$ on the subspace ${\fh}_s^{-\sigma}$.
Thus, $\nu$ is a regular special dominant weight by Lemma 6.1 in
 \cite{dp:compactification}. Moreover, the other weight vectors
appearing in $v_\sigma$ have weights $\psi$ such that $\overline{\psi}$
is of the form $\nu - \sum_{n_\alpha \ge 0,
\alpha\in\Phiplushs}^{} n_\alpha\alpha$.
Thus, by Lemma \ref{lem_dplemma}, $\overline{G\cdot v_\sigma}\cong X_\sigma.$

Note that using the Plucker embedding
of $\CGr \hookrightarrow \Proj (V),$ we can
 identify $\overline{G \cdot v_\sigma}$ with
$\overline{G_\Delta \cdot \flsigma}.$ Thus,
$\overline{G_\Delta \cdot \flsigma} \cong X_\sigma.$\qed

\begin{rem}
\label{rem_wonderful}
{\em
Let $d$ be the automorphism of the Dynkin diagram of $\fg$ such that
$\sigma = \gamma_d \Ad_{g_0}$ for some $g_0$. Consider the embedding
\[
G/G^\sigma \lrw G: \, \, gG^\sigma \Map \gamma_{d}^{-1} (g) g_0 g^{-1},
\]
which in turn gives an embedding of $G/G^\sigma$ into
the De Concini-Procesi compactification $Z_d$ of $G$.
Proposition \ref{prop_orbitequalsdp} then says that  the
closure of $G/G^\sigma$ in $Z_d$ is isomorphic to the
De Concini-Procesi compactification of $G/G^\sigma$.
}
\end{rem}

\section{The Poisson structure $\Pi_0$ on $\lag$}
\label{sec_Pi0}

\subsection{Lagrangian splittings of $\fg \oplus \fg$}
\label{sec_llgg}

By a {\it Lagrangian splitting} of $\fg \oplus \fg$ we mean a decomposition
$\fg \oplus \fg = \fl_1 + \fl_2$, where $\fl_1$ and $\fl_2$ are
Lagrangian subalgebras of $\fg \oplus \fg$. By \cite{e-l:reallag}, every
Lagrangian splitting $\fg \oplus \fg = \fl_1 + \fl_2$ gives rise to a Poisson structure
$\Pi_{{\mathfrak l}_1, {\mathfrak l}_2}$ on $\lag$ as follows: 
let $\{x_j\}$ be a basis for $\fl_1$ and 
$\{\xi_j\}$ the basis for $\fl_2$ such that
 $\la x_j, \xi_k\ra = \delta_{jk}$ for $1 \leq j, k \leq n = \dim \fg$. Set
\begin{equation}
\label{eq_Rmatrix}
R = \frac{1}{2} \sum_{j=1}^{n} (\xi_j \wedge x_j) \in \wedge^2
(\fg \oplus \fg).
\end{equation}
The action of $G \times G$ on $\lag$ defines a Lie algebra anti-homomorphism
$\kappa$ from $\fg \oplus \fg$ to the space of vector
fields on $\lag$. 
Set
\[
\Pi_{{\mathfrak l}_1, {\mathfrak l}_2} = (\kappa \wedge \kappa)(R)
=\frac{1}{2} \sum_{j=1}^{n} (\kappa(\xi_j) \wedge \kappa(x_j)).
\]

\begin{prop}
\label{prop_Pill} \cite{e-l:reallag}
For any Lagrangian splitting $\fd = \fl_1 + \fl_2$, the
bi-vector field $\Pi_{{\mathfrak l}_1, {\mathfrak l}_2}$ on $\lag$
is 
Poisson
with the property that 
 all $L_1$ and $L_2$-orbits in $\lag$ are Poisson
submanifolds with respect to $\Pi_{{\mathfrak l}_1, {\mathfrak l}_2}$, where,
for $i = 1, 2$, $L_i$ is the connected subgroup of $G \times G$ with Lie algebra
$\fl_i$.
\end{prop}

The rank of $\Pi_{{\mathfrak l}_1, {\mathfrak l}_2}$
can be computed as in the following Lemma \ref{lem_rank-Pill}. 
A version of Lemma \ref{lem_rank-Pill}  first appeared in
\cite{e-l:reallag}, and a 
generalization of Lemma \ref{lem_rank-Pill} can be found in \cite{lu-milen2}.

\begin{lem}
\label{lem_rank-Pill}
For $\fl \in \lag$, let
$\fn_{{\fg \oplus \fg}}(\fl)$ be the normalizer subalgebra of $\fl$ in $
\fg \oplus \fg$, and let  
$(\fn_{{\fg \oplus \fg}} (\fl))^\perp = \{x \in \fg \oplus \fg: 
\la x, y \ra = 0 \; \forall y \in \fn_{{\fg \oplus \fg}}(\fl)\}$. Set
\[
{{\mathcal T}}(\fl) = \fl_1 \cap \fn_{{\fg \oplus \fg}} (\fl) + 
(\fn_{{\fg \oplus \fg}} (\fl))^\perp \subset \fg \oplus \fg.
\]
Then ${{\mathcal T}}(\fl) \in \lag$, 
and the rank of
$\Pi_{{\mathfrak l}_1, {\mathfrak l}_2}$ at $\fl$ is equal to 
$\dim ({L_1}\cdot \fl) - \dim (\fl_2 \cap {{\mathcal T}}(\fl))$, where
${L_1}\cdot \fl$ is the orbit in $\lag$ of $L_1$ through $\fl$.
\end{lem}

\begin{exam}
\label{rem_D-orbits}
{\rm It is also clear from the definition of $\Pi_{{\mathfrak l}_1, {\mathfrak l}_2}$
that $\Pi_{{\mathfrak l}_1, {\mathfrak l}_2}$ is tangent to every 
$(G \times G)$-orbit in $\lag$. Thus every $(G \times G)$-orbit  in $\lag$ is a Poisson
submanifold of $(\lag, \Pi_{{\mathfrak l}_1, {\mathfrak l}_2})$,
and its closure is 
a Poisson subvariety. For example,
let $d$ be a diagram automorphism and consider the embedding of $G$
into $\lag$ as the $(G \times G)$-orbit
through $\fl_{\gamma_d}= \{(x, \gamma_d(x)): x \in \fg\}$:
\begin{equation}
\label{eq_G-into-lag}
G \hookrightarrow \lag: \; \; g \Map \{(x, \gamma_d \Ad_g(x)): x \in \fg\}.
\end{equation}
By Proposition \ref{prop_Pill}, every
Lagrangian splitting of $\fg \oplus \fg$ gives rise to a Poisson structure 
$\Pi_{{\mathfrak l}_1, {\mathfrak l}_2}$
on $G$ which  extends to  the closure $Z_d(G)$ of $G$ in $\lag$.}
\end{exam}

\begin{exam}
\label{exam_Xsigma}
{\em
By a {\it Belavin-Drinfeld splitting} we mean a Lagrangian
splitting $\fg \oplus \fg = \fl_1 + \fl_2$ in which $\fl_1 = \fg_\Delta$. 
By Corollary \ref{cor_BD}, $\fl_2$ is conjugate by an element in $G_\Delta$
to an $\fl_{S, T, d, V}$, 
where $(S, T, d, V)$ is a
Belavin-Drinfeld system (Definition \ref{dfn_BD-system}). We will also denote
a Belavin-Drinfeld splitting by 
$\fg \oplus \fg = \fg_\Delta + \fl_{BD}$, and denote by $\Pi_{BD}$ the
corresponding Poisson structure on $\lag$. By Proposition \ref{prop_Pill},
all the $G_\Delta$-orbits in $\lag$ as well as their closures
are Poisson submanifolds with
respect to any $\Pi_{BD}$. For example,
for a diagram automorphism $d$, equip $G$ and $Z_d(G)$ with
the Poisson structure $\Pi_{BD}$ via the embedding 
(\ref{eq_G-into-lag}). Then every
$d$-twisted conjugacy class in $G$, as well as its closure in
$Z_d(G)$, is a Poisson subvariety  with respect to 
every $\Pi_{BD}$. As a special case, every complex
symmetric space $G/G^\sigma$, as well as its De Concini-Procesi
compactification, inherits the Poisson structure $\Pi_{BD}$ this way.
See Remark \ref{rem_wonderful}. 
 }
\end{exam}

\smallskip
Lagrangian splittings of $\fg \oplus \fg$ up to conjugation 
by elements in $G \times G$ have been classified by P. Delorme
\cite{delorme:manin-triples}.  
A study of 
the Poisson structures $\Pi_{{\mathfrak l}_1, {\mathfrak l}_2}$ 
defined by  arbitrary Lagrangian splittings $\fg \oplus \fg
=\fl_1 + \fl_2$
will be carried out in \cite{lu-milen2}. 

For the rest of this section, we will only be concerned with the {\it standard Lagrangian splitting} of    $\fg \oplus \fg$, namely, the splitting
$\fg \oplus \fg=\fg_{\Delta}
+ \fg_{{\rm st}}^{*}$,
where 
\[
\fg_{{\rm st}}^{*}=\fh_{-\Delta} + (\fn \oplus \fn^-).
\]
We will denote by $\Pi_0$ the Poisson structure on 
$\lag$ determined by the standard Lagrangian splitting.
We will compute 
the rank of $\Pi_0$ everywhere on $\lag$. As a consequence, we will 
see that every non-empty intersection 
${{\mathcal O}} \cap {{\mathcal O}}^\prime$ of a 
$G_\Delta$-orbit ${{\mathcal O}}$
and a
$(B \times B^-)$-orbit ${{\mathcal O}}^\prime$ is a regular Poisson subvariety of $\Pi_0$, and that
the subgroup $H_\Delta = \{(h, h): h \in H\}$ of $G_\Delta$ acts transitively
on the set of all symplectic leaves in ${{\mathcal O}} \cap {{\mathcal O}}^\prime$,
where $H = B \cap B^-$.

\subsection{The rank of the Poisson structure $\Pi_0$}
\label{sec_rank-Pi0}

Let 
${{\mathcal O}}$ be a $G_\Delta$-orbit in $\lag$ and  ${{\mathcal
O}}^\prime$   a $(B \times B^-)$-orbit in $\lag$ such that
${{\mathcal O}} \cap {{\mathcal O}}^\prime \neq \emptyset$. Since
 $(\fb \oplus \fb^-) + \fg_\Delta = \fg \oplus \fg$,
${{\mathcal O}}$ and ${{\mathcal O}}^\prime$ intersect transversally
in their   $(G \times G)$-orbit.  
Since both ${{\mathcal O}}$ and
${{\mathcal O}}^\prime$ are Poisson submanifolds for $\Pi_0$, the
intersection ${{\mathcal O}} \cap {{\mathcal O}}^\prime$ is a
Poisson submanifold of $(\lag, \Pi_0)$. Thus, it is enough to
compute the rank of $\Pi_0$ as a Poisson structure in 
${{\mathcal O}} \cap {{\mathcal O}}^\prime$.
By Theorem
\ref{thm_GG-orbits-1}, there exists a generalized Belavin-Drinfeld
triple $(S, T, d)$ and  $V \in \LinLagr(\fz_S \oplus \fz_T)$ such
that ${{\mathcal O}}, {{\mathcal O}}^\prime \subset (G \times G)
\cdot \fl_{S, T, d, V}$ with $\fl_{S, T, d, V}$ given in
(\ref{eq_flSTdV}).  By Corollaries \ref{cor_BB-finite-orbits} and
\ref{cor_G-orbits-L}, there exist $w \in W$, $v, v_1 \in
W^T$, and $m \in M_{S(v, d)}$ such that
\begin{equation}
\label{eq_O} {{\mathcal O}} = G_\Delta \cdot \Ad_{(m, \dot{v})}
\fl_{S, T, d, V}, \hs 
{{\mathcal O}}^\prime = (B \times B^-) \cdot
\Ad_{(\dot{w}, \dot{v}_1)} \fl_{S, T, d, V}
\end{equation}
where $\dot{w}, \dot{v}$ and $\dot{v}_1$ are  representatives of
$w, v,$ and $v_1$ in $G$ respectively. Set
 \begin{equation}
 \label{eq_XSTdv}
 X_{S, T, d, v} = \{(z, v^{-1}z): z \in \fz_{S(v, d)},
 \gamma_d(\chi_S(z)) = \chi_T(v^{-1}z)\} + V_S \subset \fh \oplus
 \fh
 \end{equation}
 with $V_S  = \{(x, \gamma_d(x)): x \in \fh_S\}.$
One shows directly that $X_{S, T, d, v}$ is a
Lagrangian subspace of $\fh \oplus \fh$.

\begin{thm}
\label{thm_rank-Pi0}  
Let ${{\mathcal O}}$ and ${{\mathcal
O}}^\prime$ be as in (\ref{eq_O}), and suppose that
${{\mathcal O}} \cap {{\mathcal O}}^\prime \neq \emptyset$. The
rank of $\Pi_0$ at every $\fl \in {{\mathcal O}} \cap {{\mathcal
O}}^\prime$ is equal to
\[
\dim ({{\mathcal O}} \cap {{\mathcal O}}^\prime) - \dim(
\fh_{-\Delta} \cap (w, v_1) X_{S, T, d, v}),
\]
where $X_{S, T, d, v}$ is given in (\ref{eq_XSTdv}). In
particular, ${{\mathcal O}} \cap {{\mathcal
O}}^\prime$ is a regular Poisson submanifold of $\Pi_0$.
\end{thm}

\noindent
{\bf Proof.}
Let $\fl= \Ad_{(g,g)} \Ad_{(m, \dot{v})} \fl_{S, T, d, V } \in
{{\mathcal O}}$, where $g \in G$, and let ${{\mathcal T}}(\fl)$
be as in  Lemma \ref{lem_rank-Pill}. Then  the 
rank of $\Pi_0$ at $\fl$ is equal to
$\dim (G_\Delta \cdot \fl) - \dim (\fg_{{\rm st}}^{*} \cap {{\mathcal T}}(\fl))$.
It follows from the definition of ${{\mathcal T}}(\fl)$ that
$
{{\mathcal T}}(\fl) = \Ad_{(g, g)} 
{{\mathcal T}}(\Ad_{(m, \dot{v})} \fl_{S, T, d, V}).
$
Let $\fr_{S, T, d}^{\prime} = (\fn_S \oplus \fn_{T}^{-}) + \{(x,
\gamma_d(x)): x \in \fg_S\}$, and let 
\begin{equation}
\label{eq_flSTdv}
 \fl_{S, T, d, v} = X_{S, T, d, v} + \fr_{S, T,
d}^{\prime}. 
\end{equation}
By
Theorem  \ref{thm_normalizers},  
\beqa
{{\mathcal T}}(\Ad_{(m, \dot{v})} \fl_{S, T, d, V })& = &
\fg_\Delta \cap \Ad_{(m, \dot{v})} \fr_{S, T, d} + 
\Ad_{(m, \dot{v})} \fr_{S, T, d}^{\prime} \\ & = & 
(\fz_{S(v, d)}^{\prime})_\Delta + \left(\fg_{S(v,d)}^{\phi} +
\psi(\fn_v)\right)_\Delta + \Ad_{(m, \dot{v})} \fr_{S, T, d}^{\prime}.
\eeqa
Since $\Ad_{(m, \dot{v})}^{-1} \left(\fg_{S(v,d)}^{\phi} +
\psi(\fn_v)\right)_\Delta \subset \fr_{S, T, d}^{\prime}$, we have
\[
{{\mathcal T}}(\Ad_{(m, \dot{v})} \fl_{S, T, d, V }) =
\Ad_{(m, \dot{v})} \left( \Ad_{(m, \dot{v})}^{-1}(\fz_{S(v, d)}^{\prime})_\Delta
+\fr_{S, T, d}^{\prime}\right) 
=  \Ad_{(m, \dot{v})}
(\fl_{S, T, d,v}).
\]
Thus the rank of $\Pi_0$ at $\fl$ is equal to
${\rm Rank}_{\Pi_0} (\fl) = \dim {{\mathcal O}} - \dim (\fg_{\rm st}^{*} \cap
\Ad_{(gm,g\dot{v})} \fl_{S, T, d, v}).$
Let
\[
\delta = \dim((\fb \oplus \fb^-) \cap \Ad_{(gm,g\dot{v})} \fr_{S,
T, d})-\dim(\fg_{\rm st}^{*} \cap \Ad_{(gm,g\dot{v})} \fl_{S, T,
d, v}).
\]
Then
${\rm Rank}_{\Pi_0} (\fl) = \dim {{\mathcal O}}  + \delta -
\dim((\fb \oplus \fb^-) \cap \Ad_{(gm,g\dot{v})} \fr_{S,
T, d})$. 
Since
\[
\dim {{\mathcal O}}^\prime = \dim (\fb \oplus \fb^-) - \dim((\fb
\oplus \fb^-) \cap \Ad_{(gm,g\dot{v})} \fr_{S, T, d}),
\]
we have
\[
{\rm Rank}_{\Pi_0} (\fl) =\dim {{\mathcal O}} + \dim {{\mathcal O}}^\prime
+ \delta - \dim (\fb \oplus \fb^-) =\dim {{\mathcal O}} + \dim {{\mathcal O}}^\prime
+ \delta - 2\dim \fb.
\]
Since ${{\mathcal O}}$ and ${{\mathcal O}}^\prime$ intersect
transversally at $\fl$ inside the $(G \times G)$-orbit through
$\fl$, and since 
$\dim (G\times G) \cdot \fl = \dim \fg - \dim \fz_S$ 
by Proposition \ref{prop_GG-orbit-types}, we have
\beqa
{\rm Rank}_{\Pi_0}(\fl)  & = &\dim({{\mathcal
O}} \cap {{\mathcal O}}^\prime) + 
\dim ((G \times G) \cdot \fl) + \delta - 2\dim \fb \\
&= &
  \dim({{\mathcal
O}} \cap {{\mathcal O}}^\prime) - (\dim \fz_S + \dim \fh) +
\delta.
\eeqa
 It remains to compute
$\delta$. Since $\fl \in {{\mathcal O}} \cap {{\mathcal
O}}^\prime$, there exist $r \in R_{S, T, d}$ 
and $(b, b^-) \in B \times B^-$ such that $(gm,
g\dot{v}) = (b, b^-)(\dot{w}, \dot{v}_1)r$. Using
$\Ad_{(b, b^-)} (\fb \oplus \fb^-) = \fb \oplus \fb^-$ and $\Ad_{(b, b^-)}
\fg_{{\rm st}^{*}} = \fg_{{\rm st}^{*}}$, we have
\[
\delta =\dim ((\fb \oplus \fb^-) \cap \Ad_{(\dot{w}, \dot{v}_1)}
\fr_{S, T, d}) - \dim(\fg_{{\rm st}}^* \cap \Ad_{(\dot{w},
\dot{v}_1)} \fl_{S, T, d, v}).
 \]
 Set $Y = (\fn \oplus \fn^-) \cap\Ad_{(\dot{w}, \dot{v}_1)}\left(
 (\fn_{S} \oplus \fn_{T}^{-})+ {\rm span}_{{\mathbb C}} \{
 (E_\alpha, \gamma_d(E_\alpha)): \alpha \in [S]\}\right).$
Then
 \[
 (\fb \oplus \fb^-) \cap \Ad_{(\dot{w}, \dot{v}_1)}
\fr_{S, T, d} = (w, v_1) (\fz_S \oplus \fz_T + V_S) + Y.
\]
 Since $Y \subset
\fg_{{\rm st}}^{*} \cap\Ad_{(\dot{w}, \dot{v}_1)} \fl_{S, T, d,
v}$, we have $\fg_{{\rm st}}^* \cap \Ad_{(\dot{w}, \dot{v}_1)} \fl_{S, T, d, v}
= Y + \fh_{-\Delta} \cap (w, v_1) X_{S, T, d, v}.$ Thus
\beqa
\delta & = & \dim (\fz_S \oplus \fz_T + V_S) - \dim (
\fh_{-\Delta} \cap (w, v_1) X_{S, T, d, v}) \\ & = & \dim \fz_S +
\dim \fh -\dim ( \fh_{-\Delta} \cap (w, v_1) X_{S, T, d, v}).
\eeqa
 Thus    the rank of $\Pi_0$
at $\fl$ is equal to
$\dim ({{\mathcal O}} \cap {{\mathcal O}}^\prime) - \dim (
\fh_{-\Delta} \cap (w, v_1) X_{S, T, d, v})$.
\qed

\begin{rem}
\label{rem_constant-rank}
{\em
Our conclusion that ${{\mathcal O}} \cap {{\mathcal O}}^\prime$ is a
regular Poisson manifold for $\Pi_0$ follows immediately from
our computation of the rank of $\Pi_0$. It will be shown in 
\cite{lu-milen2} that a similar result holds for the Poisson structure
$\Pi_{{\mathfrak l}_1, {\mathfrak l}_2}$ on $\lag$
defined by any Lagrangian splitting $\fg \oplus \fg = \fl_1 + \fl_2$.
}
\end{rem}

\begin{cor}
\label{cor_conjclassrank}
 Equip 
$G$ with the Poisson structure $\Pi_0$ via the embedding of $G$ into 
$\lag$ in (\ref{eq_G-into-lag}) for $d = 1$. Let
$C$ be a  conjugacy class in $G$ 
and let $w \in W$ be such that
$C \cap (B^- wB) \neq \emptyset.$ Then the rank of $\Pi_0$ at every point
in $C \cap (B^- w B)$ is  
\[
\dim C - l(w) - \dim(\fh^{-w}), 
\]
where
$l(w)$ is the length of $w$, and 
$\fh^{-w}=\{x \in \fh:  w(x) = -x\}$.
In particular, $C\cap B^-B$ is an open dense leaf for $C$,
and  $\Pi_0$ is degenerate on the complement of $B^-B \cap C$ in $C$.
\end{cor}

\noindent
{\bf Proof.} By Proposition 
\ref{prop_Pill}, $C$ is a Poisson submanifold of
$(G, \Pi_0)$. By the Bruhat decomposition, 
$C = \cup_{w \in W} (C \cap (B^-wB))$.
Since $B^-B$ is open in $G$ and $C\cap B \not= \emptyset$
(Theorem 1 on P. 69 of \cite{steinberg}), 
$C\cap B^-B$ is open and dense in $C$. The rank formula
follows from Theorem \ref{thm_rank-Pi0}, and it follows easily
that $C \cap B^-B$ is a symplectic leaf and $\Pi_0$ is degenerate
on $C\cap (B^-wB)$ if $w\not= 1$.
\qed

\begin{rem}
\label{rem_nilporbitclosure}
{\em By Corollary \ref{cor_conjclassrank}, 
any unipotent conjugacy class (and its closure
in $Z_1(G)$) has an
induced Poisson structure $\Pi_0$ with an open symplectic leaf, although
the structure is not symplectic unless the orbit is a single point.
Since the unipotent variety is isomorphic to the nilpotent cone
in $\fg^*$, it follows that every nilpotent orbit in $\fg^*$ has an induced
Poisson structure with the same properties. It would be quite
interesting to compare this structure with the Kirillov-Kostant
symplectic structure.}
\end{rem}

\begin{exam}
\label{exam_double-Bruhat}
{\em
Consider the closed $(G \times G)$-orbit through a Lagrangian subalgebra of the form
$V + (\fn \oplus \fn^-)$, where $V$ is any Lagrangian subspace of $\fh \oplus \fh$.
Such an orbit can be identified with $G/B \times G/B^-$, so we can regard 
$\Pi_0$ as a Poisson structure on $G/B \times G/B^-$. Let ${{\mathcal O}}$
be a $G_\Delta$-orbit and let ${{\mathcal O}}^\prime$ be a 
$(B \times B^-)$-orbit in $G/B \times G/B^-$ such that 
${{\mathcal O}} \cap {{\mathcal O}}^\prime \neq \emptyset$. By the Bruhat 
decomposition of $G$, there are elements $w, u, v \in W$ such that
\[
{{\mathcal O}} = G_\Delta \cdot (B, wB^-), \hs 
{{\mathcal O}}^\prime = (B \times B^-) \cdot (uB, vB^-).
\]
The stabilizer subgroup of $G_\Delta \cong G$ 
at the point $(B, \dot{w}B^-) \in G/B \times G/B^-$ is $B \cap w(B^-)$.
Identify ${{\mathcal O}} \; \cong \; G/(B \cap w(B^-))$, and 
let $p: G \to {{\mathcal O}} \cong G/(B \cap w(B^-))$ 
be the projection.  It is then easy to see that 
$ {{\mathcal O}} \cap {{\mathcal O}}^\prime = p(G^{u, v}_{w}) \subset
{{\mathcal O}}$, where  
\[
G^{u, v}_{w} = (BuB) \cap (B^-vB^-w^{-1}).
\]
We will refer to $G^{u, v}_{w}$ as the {\it shifted double
Bruhat cell} in $G$ determined by $u, v$ and $w$. 
Note that $B \cap w(B^-)$ acts freely on $G^{u, v}_{w}$ by 
right multiplications, so
\[
{{\mathcal O}} \cap {{\mathcal O}}^\prime \; \cong \; 
G^{u,v}_{w} /(B \cap w(B^-)).
\]
Since $\dim {{\mathcal O}} = \dim \fg - \dim \fh - l(w)$
and $\dim {{\mathcal O}}^\prime = l(u) + l(v)$, we have
\[
\dim ({{\mathcal O}} \cap {{\mathcal O}}^\prime) = \dim {{\mathcal O}} +
\dim {{\mathcal O}}^\prime - \dim(G/B \times G/B^-) = l(u) + l(v) - l(w),
\]
and $\dim G^{u,v}_{w} = l(u) + l(v) + \dim \fh.$
By Theorem \ref{thm_rank-Pi0}, the rank of $\Pi_0$ at every point of 
${{\mathcal O}} \cap {{\mathcal O}}^\prime$ is 
\[
l(u) + l(v) - l(w) - \dim \fh^{-u^{-1}vw^{-1}},
\]
 where
$\fh^{-u^{-1}vw^{-1}} = \{x \in \fh: u^{-1}vw^{-1} x = -x\}$.
When $w = 1$, we have ${{\mathcal O}} \cong G/H$, and 
${{\mathcal O}} \cap {{\mathcal O}}^\prime \cong G^{u,v}/H$, where 
$G^{u, v} = G^{u,v}_{1}$ is the  double Bruhat cell in $G$ determined by $u$ and $v$.
The set $G^{u,v}/H$ is called a {\it reduced double Bruhat cell} in
\cite{zelevinsky:pi0}. In \cite{kogan-zelevinsky}, Kogan and Zelevinsky 
constructed toric charts on symplectic leaves of $\Pi_0$ in 
$ {{\mathcal O}} \cap {{\mathcal O}}^\prime $  (for the case when $w = 1$)
by using the so-called 
{\it twisted minors} that are developed in \cite{fomin-zelevinsky:double},
and they also constructed integrable systems on the symplectic leaves.
 It would very interesting to generalize the Kogan-Zelevinsky construction
to all symplectic leaves of $\Pi_0$ in $G/B \times G/B^-$.
  }
\end{exam}

\subsection{The action of $H_\Delta$ on the set of symplectic leaves of the
Poisson structure $\Pi_0$}
\label{sec_leaves-Pi0}

\begin{prop}
\label{prop_generalconnected-intersection}
Let $D$ be a connected complex algebraic group with connected 
algebraic subgroups
$A$ and $C$. Suppose there exists a connected algebraic subgroup
 $C_1 \subset C$
such that the multiplication morphism $A\times C_1 \to D$ is an
isomorphism to a connected open set $U$ of $D$. Let $X$ be a homogeneous space
for $D$ such that the stabilizer in $D$ of a point in $X$ is connected.
Then any nonempty intersection of an $A$-orbit in $X$ with a $C$-orbit in
$X$ is smooth and connected.
\end{prop}

\newcommand{\orbitintersect}{A \cdot x \cap C\cdot x}

\noindent
{\bf Proof.} 
Let $A \cdot x \cap C\cdot x$ be a nonempty intersection of
orbits in $X$, and note that this intersection is smooth since
the hypotheses imply that the orbits intersect transversely. 
 We show there is a fiber 
bundle $\pi:V\to U$, with fiber 
 ${\pi}^{-1}(e)\cong \orbitintersect$ over the identity and 
$V$ connected, that is trivial in the Zariski topology.
This implies the connectedness of the intersection, and hence
the proposition. The proof is inspired by
the proof of Kleiman's transversality theorem.

Let $Y=C\cdot x$ and  $Z=A \cdot x$.
 Let $h:D\times Y \to X$ be the action map and
let $i:Z\to X$ be the obvious embedding. Let 
$W = (D\times Y) {\times}_X Z$
be the fiber product. Then $h$ is a smooth fiber
 bundle (see the proof of
10.8 in \cite{ha:ag}) and the fibers $h^{-1}(x)$
are connected.
For the second claim, note that $h^{-1}(x)
= \{ (d,c\cdot x)  : dc\cdot x = x \}$ and
$\psi:h^{-1}(x) \to D_x\cdot C$ given by
$\psi(d,c\cdot x)=d$ is an isomorphism. Since
$D_x$ and $C$ are connected, the claim follows.
 Thus, the induced morphism from
$W\to Z$ also has connected fibers. Since $Z$ is connected,
it follows that $W$ is connected. Moreover, $W$ is smooth
(again by the proof of 10.8 in \cite{ha:ag}), so $W$ is irreducible.

Let $\pi:W\to D\times Y \to D$ be the composition of the induced
fiber product map with projection to the first factor.
Since $\pi^{-1}(U)$ is open in $W$, it is smooth and
irreducible, and thus connected.
Note also that $\pi^{-1}(e)\cong Y\cap Z$.
It remains to show that $\pi:\pi^{-1}(U)\to U$ is a trivial
fiber bundle. We define a free left $A$ action and a free right
$C_1$ action on $W$ by the formulas
\beqa
& a\cdot (d,y,z)=(ad, y, a \cdot z)\\ 
& c\cdot (d,y,z)=(dc, c^{-1}\cdot y, z)\\ 
& a\in A, c\in C, d\in D, y\in Y, z\in Z
\eeqa
$A$ and $C_1$ have the obvious free left and right
multiplication actions on $U$, and $\pi:\pi^{-1}(U)\to U$
is equivariant for these actions. 
%Since $U$ is a homogeneous
%space for $A \times C_1$, $\pi:\pi^{-1}(U)
%\to U$ is a fiber bundle. 
It follows that the morphism
\[
\phi: A \times C_1 \times (\orbitintersect)
\to \pi^{-1}(U), \ (a,c,v)\mapsto (ac,c^{-1}\cdot v, a\cdot v),
\ a\in A, c\in C_1, v\in Y\cap Z
\]
is a bijection, and hence is an isomorphism since
 $\pi^{-1}(U)$ is smooth.
Thus the fiber bundle is trivial. 
\qed

\begin{rem}
\label{rem_connectedintersection}
{\em
We thank Michel Brion for suggesting this approach. Note
that the Proposition \ref{prop_generalconnected-intersection}
is false as stated if 
we only assume that $A \cdot C$ is open in $D$.
For example, let $A = G_\Delta$, and let 
$C=\{(nh, h^{-1}n^-): n \in N, h \in H, n^- \in N^-\}$
be the connected subgroup
of $D=G\times G$ corresponding to $\fg_{\rm st}^*$. Let
$X = D$ and let $D$ act on $X$ by left translation.
Then the intersection of the $A$-orbit and the $C$-orbit
through the identity element of $D$ is $A \cap C$ which is
disconnected.}
\end{rem}

\begin{prop}
\label{prop_intersection-connected}
The intersection of any $G_\Delta$-orbit  and any
$(B \times B^-)$-orbit in $\lag$ is either empty or
a smooth connected subvariety of $\lag$.
\end{prop}

\noindent
{\bf Proof.}
This is a consequence of Proposition \ref{prop_generalconnected-intersection}.
Indeed, we take $A=G_\Delta$, $C = B \times B^-$, $D=G\times G$, and
$C_1 = B \times N^-$. The fact that the stabilizer of a point
in $\lag$ is connected follows from Lemma \ref{lem_RSTd-connected}.
\qed

Let $H = B \cap B^-$ 
and let $H_{\Delta} = \{(h, h): h \in H\}$. For every $G_\Delta$-orbit
${{\mathcal O}}$ and every $(B \times B^-)$-orbit ${{\mathcal O}}^\prime$
such that ${{\mathcal O}} \cap {{\mathcal O}}^\prime \neq \emptyset$,
$H_{\Delta}$ clearly leaves ${{\mathcal O}} \cap {{\mathcal O}}^\prime$ invariant.
It is easy to show that the element $R \in \wedge^2 (\fg \oplus \fg)$ 
given in (\ref{eq_Rmatrix})
is invariant under $\Ad_{(h, h)}$ for every $h \in H$. Thus the Poisson structure
$\Pi_0$ on $\lag$ is $H_\Delta$-invariant. In particular,
for every $h \in H$, $\Ad_{(h, h)} {{\mathcal E}}$ is a symplectic leaf
of $\Pi_0$ in ${{\mathcal O}} \cap {{\mathcal O}}^\prime$ if
${{\mathcal E}}$ is.

\begin{lem}
\label{lem_sigma-submersion}
Let ${{\mathcal O}}$ be a $G_\Delta$-orbit and ${{\mathcal O}}^\prime$ a
$(B \times B^-)$-orbit in $\lag$ such that
${{\mathcal O}} \cap {{\mathcal O}}^\prime \neq \emptyset$. Let
${{\mathcal E}}$ be any symplectic leaf of $\Pi_0$ in 
${{\mathcal O}} \cap {{\mathcal O}}^\prime$. Then the map
\[
\sigma: \, H \times {{\mathcal E}}  \lrw
 {{\mathcal O}} \cap {{\mathcal O}}^\prime: \, \, (h, \fl ) \Map \Ad_{(h, h)} \fl
 \]
 is a submersion.
 \end{lem}

 \noindent
 {\bf Proof.}
 Let $e$ be the identity element of $H$ and
let $\fl \in {{\mathcal E}}$. It is enough to show that
\[
\dim \ker \sigma_*(e, \fl) = \dim \fh + \dim {{\mathcal E}}_{{\mathfrak l}}
-\dim {{\mathcal O}} \cap {{\mathcal O}}^\prime,
\]
where
 $\sigma_*(e,\fl): \fh \times
T_{{\mathfrak l}} {{\mathcal E}}  \to
T_{{\mathfrak l}} ({{\mathcal O}} \cap {{\mathcal O}}^\prime)$ is the
differential of $\sigma$ at $(e, \fl)$.

We may assume that
${{\mathcal O}}$ and
${{\mathcal O}}^\prime$ are given in (\ref{eq_O}), and that
\[
\fl = \Ad_{(gm, g\dot{v})} \fl_{S, T, d, V} = \Ad_{(b\dot{w}, b^-\dot{v}_1)}
\fl_{S, T, d, V}
\]
for some $g \in G$ and $(b, b^-) \in B \times B^-$. By Theorem \ref{thm_rank-Pi0},
it is enough to show that
\[
\dim (\ker \sigma_*(e, \fl)) = 
\dim \fh - \dim(\fh_{-\Delta} \cap (w, v_1)X_{S, T, d, v}),
\]
where $X_{S, T, d, v}$ is given in (\ref{eq_XSTdv}).
Identify the tangent space of ${{\mathcal O}}$ at $\fl$ as
\[
T_{{\mathfrak l}} {{\mathcal O}} \cong \fg_\Delta / (\fg_\Delta
\cap \Ad_{(gm, g\dot{v})} \fr_{S, T, d})
\]
and let $q: \fg_\Delta \to \fg_\Delta / (\fg_\Delta
\cap \Ad_{(gm, g\dot{v})} \fr_{S, T, d})$ be the projection. Let
$p: \fg \oplus \fg \to \fg_\Delta$ be the
projection with respect to the decomposition
$\fg \oplus \fg = \fg_\Delta + \fg_{{\rm st}}^{*}$. By
the computation of ${{\mathcal T}}(\fl)$ in
the proof of Theorem \ref{thm_rank-Pi0}, the tangent
space of ${{\mathcal E}} $ at $\fl$ is given by
\[
T_{{\mathfrak l}} {{\mathcal E}} = (q \circ p) \left( \Ad_{(gm, g\dot{v})} \fl_{S, T, d, v}\right),
\]
where $\fl_{S, T, d,v}$ is given in (\ref{eq_flSTdv}).
For $x \in \fh$, let $\kappa_x$ be the vector field on
${{\mathcal O}} \cap {{\mathcal O}}^\prime$ that generates the
action of $\Ad_{(\exp tx, \exp tx)}$.
Then  $\ker \sigma_*(e, \fl) \cong \{x \in \fh:  \kappa_x(\fl) \in
T_{{\mathfrak l}} {{\mathcal E}}\}.$
Let $x \in \fh$. If $\kappa_x(\fl) \in
T_{{\mathfrak l}} {{\mathcal E}}$, then there exists
$y \in \fg$ and $(y_1, y_2) \in \fg_{{\rm st}}^*$ with
$(y+y_1, y+y_2) \in \Ad_{(gm, g\dot{v})} \fl_{S, T, d, v}$ such that
\[
(x-y, x-y) \in \fg_\Delta
\cap \Ad_{(gm, g\dot{v})} \fr_{S, T, d} = \fg_\Delta \cap
\Ad_{(gm, g\dot{v})} \fl_{S, T, d, v}.
\]
It follows that
$(x+y_1, x+y_2) \in (\fb \oplus \fb^-) \cap \Ad_{(gm, g\dot{v})} \fl_{S, T, d, v}.$
Let $r \in R_{S, T, d}$ be such that $(gm, g\dot{v}) = (b\dot{w}, b^-\dot{v}_1)r$.
Then
\[
(\fb \oplus \fb^-) \cap \Ad_{(gm, g\dot{v})} \fl_{S, T, d, v} = 
\Ad_{(b,b^-)}\left((\fb \oplus \fb^-) 
\cap \Ad_{(\dot{w}, \dot{v}_1)} \fl_{S, T, d, v}
\right).
\]
Thus there exists $(y_{1}^{\prime}, y_{2}^{\prime}) \in \fg_{{\rm st}}^*$ such that
\[
(x + y_{1}^{\prime}, x + y_{2}^{\prime}) \in
(\fb \oplus \fb^-) \cap \Ad_{(\dot{w}, \dot{v}_1)} \fl_{S, T, d, v}.
\]
If $(y^\prime, -y^\prime)$ is the $\fh_{-\Delta}$-component of
$(y_{1}^{\prime}, y_{2}^{\prime}) \in \fg_{{\rm st}}^*$, then
$(x+y^\prime, x-y^\prime) \in (w, v_1) X_{S, T, d, v}.$
Thus $(x, x) \in p((w, v_1) X_{S, T, d, v})$, where
$p$ denotes the projection $\fh \oplus \fh \to \fh_\Delta$
along $\fh_{-\Delta}$. 
Conversely, if $x \in \fh$ is such that
  $(x, x) \in p((w, v_1) X_{S, T, d, v})$, then
there exists $y^\prime \in \fh$ such that
\[
(x+y^\prime, x-y^\prime) \in (w, v_1) X_{S, T, d, v} \subset \Ad_{(\dot{w},
\dot{v}_1)} \fl_{S, T, d, v},
\]
and thus $\Ad_{(b, b^-)} (x+y^\prime, x-y^\prime) \in
\Ad_{(gm, g\dot{v})} \fl_{S, T, d, v}.$
Since 
\[
\Ad_{(b, b^-)} (x+y^\prime, x-y^\prime) =
 (x+y^\prime, x-y^\prime) {\rm mod} (\fn \oplus \fn^-),
 \]
we see that $p(\Ad_{(b, b^-)} (x+y^\prime, x-y^\prime)) = (x, x),$
so
$\kappa_x(\fl) \in T_{{\mathfrak l}} {{\mathcal E}}$.
Thus we have shown that
\[
\ker \sigma_*(e, \fl) \cong \{ x \in \fh: (x, x) \in p((w, v_1) X_{S, T, d, v})\}.
\]
It follows that
 $\dim (\ker \sigma_*(e, \fl)) = \dim \fh - \dim (\fh_{-\Delta} \cap
(w, v_1) X_{S, T, d, v}).$
The lemma now follows from Theorem \ref{thm_rank-Pi0}.
\qed

\begin{thm}
\label{thm_H-translates-leaves}
For every $G_\Delta$-orbit
${{\mathcal O}}$ and every $(B \times B^-)$-orbit ${{\mathcal O}}^\prime$
such that ${{\mathcal O}} \cap {{\mathcal O}}^\prime \neq \emptyset$,
 $H_\Delta$ acts transitively on the set of symplectic leaves of $\Pi_0$ in
 ${{\mathcal O}} \cap {{\mathcal O}}^\prime$.
\end{thm}

\noindent
{\bf Proof.} For $\fl \in
{{\mathcal O}} \cap {{\mathcal O}}^\prime$, let  ${{\mathcal E}}_{{\mathfrak l}}$
be the symplectic leaf of $\Pi_0$ through $\fl$, and let
\[
{{\mathcal F}}_{{\mathfrak l}} = \bigcup_{h \in H} \Ad_{(h, h)}
{{\mathcal E}}_{{\mathfrak l}} \subset
{{\mathcal O}} \cap {{\mathcal O}}^\prime.
\]
Then it is easy to see that either ${{\mathcal F}}_{{\mathfrak l}}
\cap {{\mathcal F}}_{{\mathfrak l}^\prime} = \emptyset$ or
${{\mathcal F}}_{{\mathfrak l}}
={{\mathcal F}}_{{\mathfrak l}^\prime}$ for any $\fl, \fl^\prime \in
{{\mathcal O}} \cap {{\mathcal O}}^\prime$.
 It  follows from Lemma \ref{lem_sigma-submersion} that
${{\mathcal F}}_{{\mathfrak l}}$ is open in
  ${{\mathcal O}} \cap {{\mathcal O}}^\prime $ for every $\fl$.
  Since ${{\mathcal O}} \cap {{\mathcal O}}^\prime $ is connected
by Proposition \ref{prop_intersection-connected},
    ${{\mathcal O}} \cap {{\mathcal O}}^\prime
  ={{\mathcal F}}_{{\mathfrak l}}$ for every $\fl \in
{{\mathcal O}} \cap {{\mathcal O}}^\prime $.

\qed

\section{Lagrangian subalgebras of $\fg\oplus\fh$}

Let again $\fg$ be a complex semi-simple Lie algebra with
Killing form $\llgg$. Let $\fh \subset \fg$ be a
Cartan subalgebra. In this section, we will consider
the direct sum Lie algebra $\fg\oplus\fh$, together with
the symmetric, non-degenerate, and ad-invariant bilinear form
\begin{equation}
\label{eq_lara-fgfh}
\la (x_1,y_1),(x_2, y_2)\ra=\ll x_1,x_2 \gg - \ll y_1,y_2\gg, \hs
x_1, x_2 \in \fg, y_1, y_2 \in \fh.
\end{equation}
We wish to describe the variety $\Lagr(\fg\oplus\fh)$
of Lagrangian subalgebras of $\fg\oplus\fh$
with respect to $\lara$. We can describe all such Lagrangian
subalgebras by using a theorem of
Delorme \cite{delorme:manin-triples}.

\begin{dfn} \cite{delorme:manin-triples}
\label{dfn_finvolution}{\em Let $\fm$ be a complex reductive Lie
algebra with simple factors ${\fm}_i, i\in I$.
 A complex linear involution $\sigma$ of $\fm$ is called
an $f$-involution if $\sigma$ does not preserve any ${\fm}_i$.}
\end{dfn}

\begin{thm}\cite{delorme:manin-triples}
\label{thm_delorme-1} Let $\fu$ be a complex reductive Lie
algebra with a symmetric, non-degenerate, and ad-invariant
bilinear form $\beta$.

1). Let $\fp$ be a parabolic subalgebra of $\fu$
with Levi decomposition $\fp = \fm + \fn$, and decompose
$\fm$ into
$\fm = \fmss + \fz$, where $\fmss$ is its semisimple part and
$\fz$ its center.  Let
$\sigma$ be an $f$-involution of $\fmss$ such that ${\fmss}^\sigma$ is
a Lagrangian subalgebra of $\fmss$ with respect to the
restriction of $\beta$, and let $V$ be a Lagrangian subspace of $\fz$
with respect to the restriction of $\beta$. Then 
$\fl(\fp,\sigma,V):={\fmss}^{\sigma} \oplus V \oplus \fn$ 
is a Lagrangian subalgebra
of $\fu$ with respect to $\beta$. 

2). Every Lagrangian subalgebra of 
$\fu$ is $\fl(\fp,\sigma,V)$
for some $\fp$, $\sigma$, and $V$ as in 1).
\end{thm}

\begin{prop}
\label{prop_gplushlagr}
Every Lagrangian subalgebra of $\fg\oplus\fh$ with respect to 
$\lara$ given in (\ref{eq_lara-fgfh}) is of the form
$\fn +V$, where $\fn$ is the nilradical of a Borel 
subalgebra $\fb $ of $\fg$, $V$ is a Lagrangian subspace
of $\fh \oplus \fh$, and $\fn + V =\{(x + y_1, y_2): x \in \fn, (y_1, y_2) \in V\}.$
\end{prop}

\smallskip
\noindent{\bf Proof.} Applying Delorme's theorem to our case
of $\fu = \fg \oplus \fh$ and $\lara$ as the bilinear form
$\beta$,
  every Lagrangian subalgebra of $\fg \oplus \fh$
is of the form
\[
\fl = \{(x + y_1, y_2): \, x \in \fmss^\sigma + \fn, \,
(y_1, y_2) \in V\}
\]
for some parabolic subalgebra $\fp$ of $\fg$ with Levi decomposition
$\fp = \fm + \fn = \fmss + \fz + \fn$, an $f$-involution $\sigma$
on $\fmss$,
and a Lagrangian subspace $V$ of $\fz \oplus \fh$. We will now show that 
if $\fmss \neq 0$ and if 
$\sigma$ is an $f$-involution of $\fmss,$ then ${\fmss}^\sigma$ is not
an isotropic subspace of $\fmss$ for 
the restriction of the Killing form $\llgg$ of $\fg$ to $\fmss$. It
follows that $\fp$ must be Borel, which gives Proposition
\ref{prop_gplushlagr}.

 Assume that $\fmss \neq 0$. Let $\fm_i$
be a simple factor of $\fmss$. Then since $\fm_i$ is simple, it has
a unique nondegenerate invariant form up to scalar multiplication.
Hence the Killing form $\llgg$ of $\fg$ restricts to a scalar multiple of
the Killing form of $\fm_i$. Recall that the Killing form on a
maximal compact subalgebra of a semisimple Lie algebra is negative
definite. It follows that the Killing form of $\fg$ restricts to
a nonzero positive scalar multiple of the Killing form on $\fm_i$.
Suppose that $\sigma$ is an involution of $\fmss$ mapping
$\fm_i$ to $\fm_j$ with $i \neq j$.  
Then $\sigma$ is an isometry with respect to the Killing
form of $\fm_i$ and the Killing form of $\fm_j$. Thus, there exists
a nonzero positive scalar $\mu$ such that
$\ll \sigma(x), \sigma(y) \gg = \mu \ll x, y \gg, \forall 
x, y \in \fm_i.$
The fixed point set ${\fmss}^{\sigma}$ contains the subspace
$\{ x + \sigma(x): \, x\in \fm_i \}$. Let $x$ be a nonzero
element of a maximal compact subalgebra of $\fm_i$. Then
$\ll x + \sigma(x), x + \sigma(x) \gg = (1+\mu) \ll x, x, \gg \neq 0$.
Thus $\fmss^\sigma$ cannot be isotropic with respect to $\llgg$.
\qed

Let $G$ be the adjoint group of $\fg$, and let 
$B$ be the Borel subgroup of $G$ corresponding to a
Borel subalgebra $\fb$. 
\begin{thm}
\label{thm_lagrfgfh}
The variety
$\Lagr(\fg\oplus\fh)$ is isomorphic to the trivial fiber bundle
over
$G/B$ with fibre $\LinLagr(\fh\oplus\fh,\lara)$. 
In particular, $\Lagr(\fg\oplus\fh)$ is
smooth with two disjoint irreducible components, 
corresponding to the two connected
components of $\LinLagr(\fh\oplus\fh,\lara)$.
\end{thm}

\smallskip
\noindent {\bf Proof. } Identify $G/B$ with the variety of all
Borel subalgebras of $\fg$. We map $\Lagr(\fg\oplus\fh)$ to $G/B$
by mapping a Lagrangian algebra $\fl=\fn + V$ to the unique Borel
subalgebra  
with nilradical $\fn$. The fiber over $\fn$ may be identified
with $\LinLagr(\fh\oplus\fh, \lara)$.
The claim about connected components follows from the fact the
bundle is trivial.
\qed

\end{document}